\documentclass[11pt,halfline,a4paper]{ouparticle}

\usepackage{amsmath,amssymb,amsthm}
\usepackage[numbers,sort&compress]{natbib}
\usepackage{subfig}
\usepackage{booktabs}

\usepackage{cleveref}
\usepackage{todonotes}

\usepackage{bbm}
\AtBeginDocument{%
  \renewcommand{\H}{\mathcal{H}}
  \renewcommand{\L}{\mathcal{L}}
}
\newcommand{\R}{\mathbb{R}} 
\newcommand{\E}{\mathbb{E}} 
\newcommand{\No}{\mathbb{N}_0} 

\newcommand{\LL}{\mathcal{L}}

\newcommand{\G}{\mathbb{G}} 
\renewcommand{\P}{\mathbb{P}} 
\newcommand{\PP}{\mathbbm{P}}

\newcommand{\EE}{\mathbbm{E}}

\newcommand{\cG}{\mathcal{G}}

\newcommand{\cV}{\mathcal{V}}

\newcommand{\cE}{\mathcal{E}}

\newcommand{\cM}{\mathcal{M}}

\newcommand{\bY}{\mathbf{Y}} 
\newcommand{\bZ}{\mathbf{Z}} 
\newcommand{\testzero}{\mathcal{H}_d}

\DeclareMathOperator{\Var}{Var}

\newcommand{\tX}{{\tilde{X}}}

\newcommand{\gr}[1]{\textcolor{magenta}{#1}}
\newcommand{\wx}[1]{\textcolor{brown}{#1}}
\newcommand{\xz}[1]{\textcolor{blue}{#1}}

\newcounter{theorem}
 \newtheorem{Theorem}[theorem]{Theorem}
 
 \newcounter{lemma}
 \newtheorem{Lemma}[lemma]{Lemma}
 
 \newcounter{corollary}
 \newtheorem{Corollary}[corollary]{Corollary}
 
 \newcounter{proposition}
 \newtheorem{Proposition}[proposition]{Proposition}
 
 \newcounter{characterization}
 \newcounter{property}
 \newcounter{problem}
 \newcounter{example}
 \newtheorem{Example}[example]{Example}
 
 \newcounter{remark}
 \newtheorem{Remark}[remark]{Remark}
 
 \newcounter{definition}
 \newtheorem{Definition}[definition]{Definition}
 
 \newcounter{hypothesis}
 \newcounter{notation}
 \newcounter{assumption}
\newcommand{\paperstatement}[2]{%
  \par\addvspace{\baselineskip}%
  \begingroup
  \emergencystretch=3em\relax
  \noindent\textbf{#1:} #2\par
  \endgroup
}

\providecommand{\authorcontributions}[1]{%
  \paperstatement{Author Contributions}{#1}%
}

\providecommand{\funding}[1]{%
  \paperstatement{Funding}{#1}%
}

\providecommand{\institutionalreview}[1]{%
  \paperstatement{Institutional Review Board Statement}{#1}%
}

\providecommand{\dataavailability}[1]{%
  \paperstatement{Data Availability Statement}{#1}%
}

\providecommand{\acknowledgments}[1]{%
  \paperstatement{Acknowledgments}{#1}%
}

\providecommand{\conflictsofinterest}[1]{%
  \paperstatement{Conflicts of Interest}{#1}%
}

\providecommand{\appendixtitles}[1]{}
\providecommand{\appendixstart}{}

\providecommand{\reftitle}[1]{%
  \renewcommand{\refname}{#1}%
}

\begin{document}

\title{Approximating 
Distributions on Vertex Degrees 
in  Random Hypergraphs via Stein's Method}



\author{%
\name{Xianghe Zhu}
\address{Department of Statistics, London School of Economics}
\email{X.Zhu58@lse.ac.uk}
\and
\name{Wenkai Xu}
\address{Department of Statistics, University of Warwick}
\email{Wenkai.Xu@warwick.ac.uk}
\and
\name{Gesine Reinert}
\address{Department of Statistics, University of Oxford}
\email{reinert@stats.ox.ac.uk}}

\abstract{Complex data are increasingly represented as hypergraphs. We establish bounds for approximating joint vertex-degree distributions by  independent multivariate Poisson and correlated normal distributions  in   hypergraph models  with independent hyperedges. We discuss the approximation of the distribution under different density regimes which are further supported by simulation studies.}

\date{October 29, 2026}

\keywords{vertex degrees; Erd\"os-R\'enyi hypergraphs; Stein's method;  multivariate Poisson 
approximation; multivariate 
normal approximation}

\maketitle

\section{Introduction}
\label{sec:intro}

Relational data and relationships between entities are commonly studied in the form of graph data (or network data). Beyond graph representations which may only capture particular relationships between two entities through vertex pairs, multiple-entity relations can be described using \textit{hypergraphs}.
As a generalisation over graph, hypergraphs can describe more complex relationships and interactions between more than two entities. In particular, the \textit{hyperedge} captures such a relationship via a vertex-tuple, allowing the order of the tuple being  greater or equal to two.
Hypergraphs are increasingly used to model complex dependencies, 
including protein complex 
interactions \citep{klimm2021hypergraphs}, cellular connections \citep{klamt2009hypergraphs} or co-authorship in scientific collaborations \citep{ji2016coauthorship}. 
In these situations, pairwise interactions between entities (vertices), as reflected by an adjacency matrix of a graph, may not capture enough information and interactions between multiple vertices can be represented as hyperedges.
Probabilistic models \citep{darling2005structure, chodrow2020configuration, mukherjee2022estimation}, statistical tools \citep{austin2008exchangeable,ghoshal2009random} and machine learning methods \citep{zhou2006learning, feng2019hypergraph, gao2020hypergraph} are developed for hypergraphs to study the relational
data of these types. 

In order to understand the typical behaviour of hypergraph summary statistics, their (approximate) distribution under simple models is of interest. In this work, we focus on the joint distribution of vertex degrees in a random hypergraph model with independent hyperedges. A difficulty arises in that even in models with independent hyperedges, the hyperedge degree counts are not independent due to the possibility of shared hyperedges.

While our probabilistic models 
assume independence between the presence of hyperedges, the departure of observed  distributions on 
hyperedge degree distributions from the desired 
ones may indicate dependency between hyperedges and researchers would further investigate more involved structures in the hypergraph.
In particular, we employ Stein's method \citep{stein1972bound} to study the distribution 
approximations
for the hypergraph where we provide 
explicit, non-asymptotic bounds, in terms of appropriate distance measures,
when comparing with
Poisson distributions \citep{barbour1992poisson} as well as 
normal distributions \citep{chen2011normal}. 
Our approximation results 
provide analysis and interpretation for hypergraph behaviours in different {hypergraph} density regimes.
Our theoretical analysis is illustrated with simulations. 

This paper is organised as follows. 
In \Cref{sec:background}, we introduce relevant notations and provide illustrations for the class of hypergraphs we focus on in this work. Then we present our main approximation results for degrees in \Cref{sec:deg}; simulations to illustrate the results are provided  
in \Cref{sec:simulation}. We conclude the paper by discussing the limitations on assumptions in this work and potential improvements in the future works in \Cref{sec:discussion}. All proofs are deferred to the appendix. The code for the experiments is available at \xz{\url{https://github.com/xianghezhu/approx_hyper}}.

\section{Preliminaries}\label{sec:background}

To measure distances in distributions we employ instances of integral probability metrics (IPMs). An {\it integral probability metric} indexed by a class $\mathcal{H}$ is of the form
$$d(p,q) = \sup_{f \in \mathcal{H}} \Big| \int f dp - \int f dq \Big| $$
where $p$ and $q$ are probability distributions. When $\H = \{h: 0\leq h(x) \leq 1, x\in \No{^M} \} $, {for some positive integer $M$}, the IPM is called {\it total variation distance,} denoted $d_{TV}$. When $\H = \{\mathbbm{1} (x\leq z), z\in \R\} $, the IPM is called {\it Kolmogorov distance,} denoted $d_{K}$; $\mathbbm{1}_A (x) =  \mathbbm{1} (x \in A)$ denotes the indicator function which is 1 if $x \in A $ and 0 otherwise. 
When $\H$ is the set of all indicator functions of convex Borel
subset of $\R^d$, the IPM is called the {\it convex set distance} $d_{{\mathcal C}}$. 
The convex set distance is a natural extension of the one-dimensional Kolmogorov distance; see for example the discussion in \citep[Section 2.1]{meliot2022gaussian}. 
We denote by $\L(W)$  the probability law of the random variable $W$. The Poisson distribution with parameter $\lambda$ is abbreviated $Po(\lambda)$.

Denote a vertex-labelled hypergraph on a vertex set $\cV = [n]:= \{1,\dots,n\}$ as an ordered pair $\cG_n = (\cV, \cE)$ where $\cE$ is a set of tuples of order at least $2$ that consists of elements from $V$, to represent the hyperedges. 
We restrict our interest to the hypergraphs that are defined
on a finite set of vertices $\cV$ with $|\cV|<\infty$; and are \textit{simple}, i.e. the hyperedges $e\in \cE$ are undirected (a hyperedge tuple $e\in \cE$ has no particular order), unweighted, {and there are} no self-loops ($\forall e \in \cE$, and $\{v_1,\dots,v_k\} \in e$, $v_i\neq v_j, \forall 1\neq i<j \leq k$) 
{and} no multiple hyperedges.
A hypergraph in which all hyperedges have exactly $k$ elements is a {\it $k$-uniform} hypergraph; $2$-uniform hypergraphs are classical graphs.

\begin{figure}[t!]
    \centering
    \includegraphics[width=0.75\textwidth]{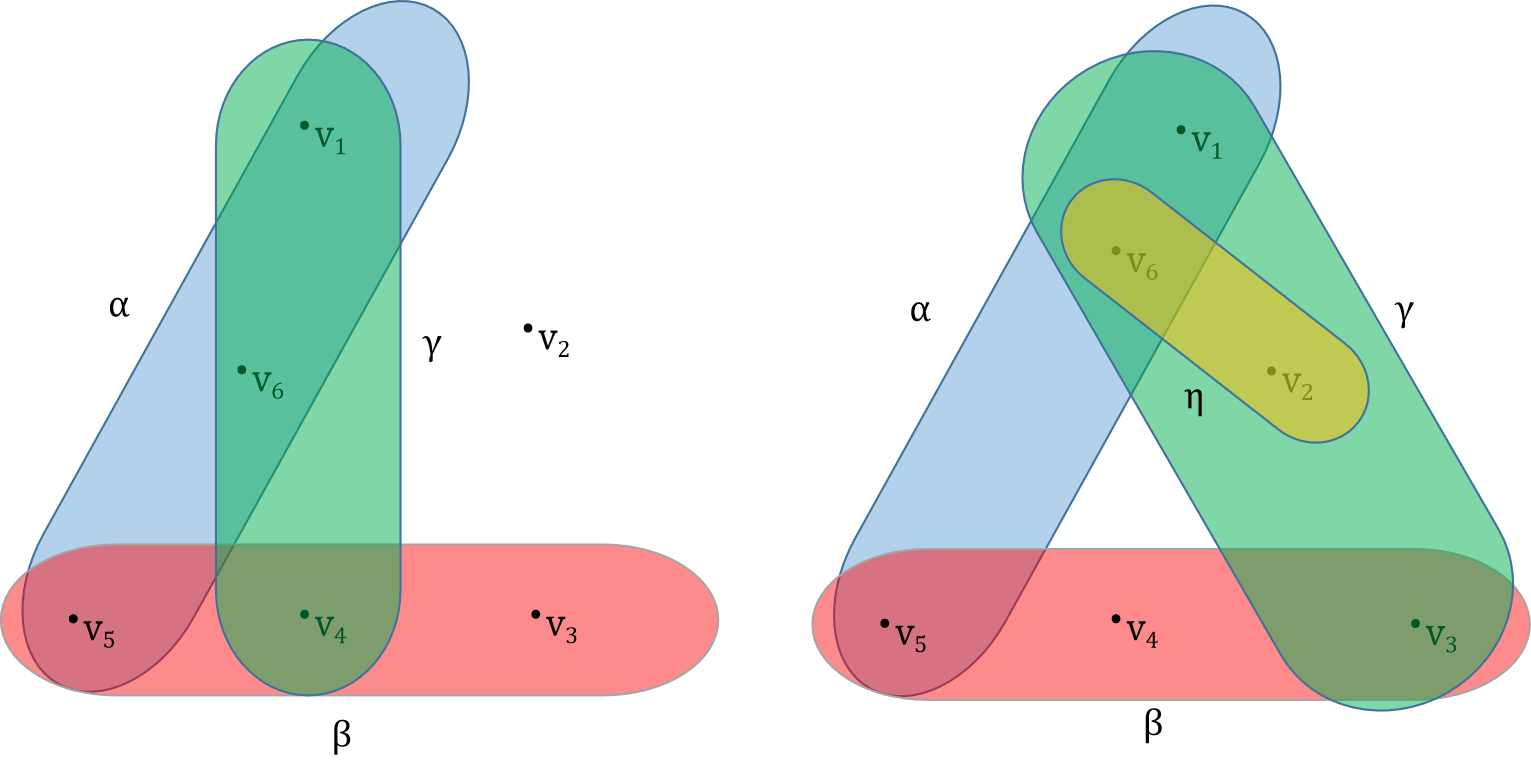}
    \caption{left: a sample from $\mathbb{G}_{6, 0.06}^{(3)}$; right: a sample from $\mathbb{G}_{6, (0.02, 0.04, 0.02)}$     \label{fig:tri}}
\end{figure}

In this work, we focus on the class of random hypergraph models where the presence of hyperedges are assumed be independent; 
graph counterparts include Bernoulli random graphs, 
Erd\"os-R\'enyi mixture models, and Chung-Lu models; see for example Chapter 4 in \cite{barbour2026networks}. 
These seemingly simple probabilistic models can embed rich statistical properties and are worth detailed investigation, which we focus in this work. 
In this paper we use the following random hypergraph model, called \emph{(non-uniform) Erd\"{o}s-R\'{e}nyi (random ) hypergraph} in \cite{schawe2022higher, agostinelli2026higher}. The same construction is also called non-uniform binomial random hypergraph in  
Cooley et al.~\cite[Definition~1.2]{cooley2020phase}.

\begin{Definition}\label{def:ER}
A  non-uniform 
Erd\"{o}s-R\'{e}nyi (random) 
hypergraph model $\mathbb{G}_{n, (p_2,\cdots, p_K)}$ of order $K{\in\{2,3,\cdots, n\}}$ is a 
hypergraph model with $|\cV| = n$ and for $k \in \{2,\cdots, K\}$, each hyperedge of order $k$ has probability $p_k$ to be present, independently of all other hyperedges. Hyperedges of order larger than $K$ have probability zero.
\end{Definition}

A $\mathbb{G}_{n, (p_2,\cdots, p_K)}$ hypergraph in which $p_k=0$ except for $k=K$ is a {\it $K$-uniform hypergraph}, denoted $\mathbb{G}^{(K)}_{n, p_K}$.  In \Cref{fig:tri}, we show Erd\"{o}s-R\'{e}nyi hypergraph samples of $6$ vertices, one from a 3-uniform hypergraph with $p_3=p=0.06$ and one non-uniform hypergraph  with $p_2 = 0.02, p_3 = 0.04, p_4 = 0.02$. 
The hyperedges are represented by the color bands and marked as Greek letters. The vertex coverage in each color band indicates the order of the hyperedge. 
In $\G^{(3)}_{6, 0.06}$, there are ${6 \choose 3}=20$ possible hyperedges of order $3$. 
With $p=0.06$, the expected number of hyperedges of order $3$ is $20 \times 0.06 = 1.2$; the depicted sample has more than the expected number hyperedges. Similarly, in $\G_{6, (0.02, 0.04, 0.02)}$, the expected number hyperedges of order $k$ is ${6\choose k}\times p_k$, that is $0.3, 0.8, 0.3$ for $k=2,3,4$, while the sample includes $1,2,1$ hyperedges respectively. 

\section{
Distributional Comparison Results for Joint Vertex Degrees}\label{sec:deg}

Let $\cE=\cE_n$ denote the set of possible simple hyperedges of size at least $2$ on $n$ vertices.
 The  {\it (hyperedge) degree} of a vertex $v$ in a hypergraph $\cG=(\cV,\cE)$ is 
\begin{align}
    \label{eq:degree}
{\rm{deg}}(v) = \sum_{e \in \cE} \mathbbm{1}( v \in e) = 
\sum_{e \in E: v\in e} \mathbbm{1}( e \in \cE) ,\end{align}
the number of hyperedges that contain $v$. 
In the Erd\"{o}s-R\'{e}nyi random hypergraph model $\G_{n, (p_2,\cdots,p_K)}$ of order $K$, the largest possible hypergraph edge size is $K$, and the degree of a vertex is a sum of $N=N_K = \sum_{k=2}^K  {n-1 \choose k-1}$  independent random variables, with mean 
\begin{align}\label{eq:mean}
\lambda = \mathbbm{E}( {\rm{deg} }(v)) = \sum_{k=2}^K {n-1 \choose k-1} p_k. 
\end{align}
Throughout this paper we assume that $\lambda >0$.
From 
\citep[Eq.(1.23)]{barbour1992poisson} it follows immediately that   in $\mathbb{G}_{n, (p_2,\dots,p_K)}$, we have \begin{equation} \label{eq:poissonbound}  d_{T V}(\mathcal{L}(deg(v)), Po(\lambda))  \leq 
\min(1,\lambda^{-1})\sum_{k=2}^K {n-1 \choose k-1} p_k^2.
\end{equation}
The Poisson bound \eqref{eq:poissonbound} reveals the degree behavior for different regimes w.r.t. hyperedge order $k$. In a $k$-uniform hypergraph, if   $p\in [0,{{n-1 \choose k-1}}^{-1}]$, {for $k$ fixed} the upper bound for the total variation  distance goes to zero when $n\to \infty$, the rate for the total variation  distance to vanish is $O(n^{-k+1})$.  
The bound is consistent with the result for Bernoulli random graph setting by taking $k = 2$ in \citep[Eq.(1.8)]{barbour1992poisson}.

Next we investigate the distributions on joint vertex degrees of vertices in a set $\cM \subset \cV$ of size $M$. 
Here a complication arises from the dependence between degrees. To make this transparent, for $v \in \cV$ denote the set of all possible hyperedges that $v$ can belong to by 
\begin{align*}
        \Gamma_v = \{ (v, e): e \in \cE, v \in e{, |e| \le K} \},
\end{align*}
which is of size $N{=N_K}$, and 
set 
\begin{align*}
 \Gamma =    \Gamma_{\cM} = \cup_{v \in \cM} \Gamma_v, 
\end{align*}
which is of size ${MN}$, noting that $\Gamma_u \cap \Gamma_v = \emptyset$ for $u \ne v$. Further, 
for $\alpha = (v, e) \in \Gamma$ set 
\begin{align*}
    X_\alpha = \mathbbm{1}( e \in \cE).  
\end{align*}
Then for $\beta = (u,f)$, $X_\alpha$ and $X_\beta$ are independent if $e \ne f$, due the independence of the hyperedge indicators, but for $\beta = (u,e)$ with $u \ne v$ and $u,v \in e$ we have 
\begin{align} \label{eq:dep}
    \E X_\alpha X_\beta =
    \P  (e \in \cE ). 
\end{align}
For $\alpha = (v,e)$ we abbreviate 
\begin{align*}
    \mu_\alpha =  \P (e \in \cE) = 
    \sum_{k=2}^{K} \mathbbm{1}(|e| =k) \,  p_k. 
\end{align*}
Moreover for $\alpha = (v,e) \in \Gamma$ we set 
\begin{align*}
    B_\alpha = \{ (u, e): (u,e) \in \Gamma\}, 
\end{align*}
of size $|B_\alpha| ={|e\cap\cM|\leq} {\min\{|e|, M\} \le M}$; 
extending the pairwise independence argument it is easy to see 
 that $X_\alpha$ is independent of $\{ X_\gamma, \gamma \not \in B_\alpha\}.$
Hence the dependence between the degrees is local: each $X_\alpha$ depends only on a small number of other hyperedge indicators. This observation is put to use in the Chen-Stein method for Poisson process approximation. Using \citep[Theorem 10.A]{barbour1992poisson} we obtain the following result.

\begin{Theorem}\label{thm:jointp}
Let $\mathbf{Y} \in \mathbb{R}^{M}$ be the vector of the joint vertex degrees in $\mathbf{G}_{n,(p_2,\dots,p_K)}$, i.e. $Y_v, v\in V$ denotes the degree of vertex $v$, where $v \in \cM$ and $
|\cM |= M$. 
Then
\begin{equation*}
  d_{T V}\left(\mathcal{L}\left(\left\{Y_v\right\}_{v=1}^M\right), \prod_{v=1}^M \operatorname{Po}(\lambda) \right) \leq {\min \{ 1; B_P\}} 
\end{equation*}
where 
\begin{equation}\label{eq:B}
     B_P:= M\sum_{k=2}^{K}  {n-1 \choose k-1}p_k^2+M(M-1) \sum_{k=2}^{K}  {{n-2} \choose {k-2}} p_k(1+p_k).
\end{equation}
\end{Theorem}

\begin{Example}\label{ex:Poisson}
\begin{enumerate}
    \item 
When the hypergraph is indeed a Bernoulli random graph, so that $K=2$ and $p_k=0$ for $k \ge 3$, and $k=2$ is the only possible edge size, then the bound can be simplified to 
$ M ((n-1) p^2 + (M-1)p(1+p))$.
Thus, the smaller the set $\cM$, the better the approximation by an independent Poisson process. 
 In particular the bound will be small for large $n$ when $p = \frac{c}{n}$ for some constant $c$ and $M\le n^\alpha$ for a constant $\alpha  < {\frac{1}{2}}$.
\item As another example, if $K=n$ and  $p_k= c n^{-k+1}$ for all $k$ then
$\lambda =  c \sum_{k=2}^n {n-1 \choose k-1} n^{-k+1} = c\left((1+\frac{1}{n})^{n-1}-1\right)\approx c (e -1)$ as $ n \rightarrow \infty$.
Then, the bound in Theorem \ref{thm:jointp} is of order $M^2n^{-1}$. Thus, the bound will remain small even if $M = o ( n^{1/2}) $ grows slowly with $n$.
\item 
If we take $K=n$ and $p_k = c / {{n-1} \choose {k-1}} $ then $\lambda = \sum_{k=2}^n {n-1 \choose k-1} p_k=c(n-1)$ and the first term of the bound in Theorem \ref{thm:jointp} behaves like $Mc^2 \sum_{j=1}^{n-1} ( {{n-1} \choose j})^{-1} $; for constant $c$ and $M$ this term is 
of order 1.  For $M\geq 2$ the second term is of order $M^2\sum_{k=2}^n{n-2\choose k-2}p_k = M^2\frac{c}{n-1}\sum_{k=2}^n(k-1)=\frac{cnM^2}{2}$. So for $M\geq 2$ the  uncapped bound $B_P$ is of order $n$. 
\item If instead we take $K=3$ and $p_k = c / {{n-1} \choose {k-1}} $ then $\lambda = 2 c $ and the bound in Theorem \ref{thm:jointp} is of order $M^2/n$. Hence this bound will be small as $n \rightarrow \infty$ as long as $M(n) = {o(n^{1/2}}).$  
\end{enumerate}
\end{Example}

\begin{Example}\label{ex:degreesequence} {Application: 
{approximating the sampled} degree sequence}

If we think of $\cM$ as a sample of size $M$ from the vertex set, drawn independently of the hypergraph, then we may want to use the sample to estimate the degree distribution in the hypergraph; this is the vector $(d_0, d_1, \ldots, d_{L})$ with
$d_\ell= \frac1n \sum_{v=1}^n \mathbbm{1}({\rm deg} (v) = \ell)$ the proportion of vertices with hyperedge degree $\ell$, $\ell=0, \ldots, L  := \sum_{\ell=2}^n {n-1 \choose \ell-1}{= 2^{n-1} -1}$.
Theorem \ref{thm:jointp} 
gives the bound in total variation by the vertex degrees of the vertices in $\cM$ replaced by i.i.d.\,Poisson variables $Z_1, \ldots, Z_M$ with mean $\lambda$. Let $q_\ell =Po(\lambda)\{ \ell\}$ be the probability that a Poisson variable with mean $\lambda$
takes on the value $\ell$. 
Letting  
$$ Q_\ell = \frac1M \sum_{v\in \cM} \mathbbm{1}(Z_v= \ell)$$
we note that $ \sum_{v\in \cM} \mathbbm{1}(Z_v= \ell)$ has the binomial distribution $Bin(M, q_\ell)$. Recall that the total variation bound is a bound on the distance in probabilities over any possible sets. Using the total variation bound we get a bound on the joint distribution of members of the degree distribution. With
$d_\ell (\cM)= \frac1M \sum_{v \in \cM} \mathbbm{1}({\rm deg} (v) = \ell)$, for every degree value $\ell \ge 0$ and  any non-negative integer $m \in \{0, 1, \ldots, M\}$,
employing the bound in 
in 
Theorem \ref{thm:jointp}, 
\begin{align*}
 \PP ( Q_\ell \le {m}{/M}) - {\min\{1,B_P\}}   \le \PP (d_\ell (\cM) \le {m}{/M} ) 
\le \PP ( Q_\ell \le {m}{/M}) + {\min\{1,B_P\}}; 
\end{align*}
the upper and lower bounds can be calculated easily using the binomial distribution. Here,  $B_P$ is given in \eqref{eq:B}.
Similarly, for any vector ${\bf L} = (\ell_1, \ldots, \ell_{{q}})$ of $q$ distinct non-negative integers, the joint distribution of $({M} Q_\ell, \ell \in {\bf L}; M - \sum_{\ell \in {\bf L}} M Q_\ell)$ is multinomial
${\rm Mult}(M; q_{\ell},  \ell \in {\bf L};   1 - \sum_{ \ell \in {\bf L}}q_{\ell})$, and for any 
non-negative integers $m_\ell, \ell \in {\bf L}$, 
\begin{align*}
  \PP ( Q_{\ell} \le m_\ell / M, {\ell \in {\bf L}} ) - {\min\{1,B_P\}}  & \le \PP (d_{\ell} (\cM) \le m_{{\ell}} / M , {\ell \in {\bf L}} 
  ) \\
  & 
\le \PP ( Q_{\ell} \le m_\ell / M,  {\ell \in {\bf L}} )
+ \min\{1,B_P\}.    
\end{align*} 
\end{Example}

\subsection{Multivariate Normal Approximation for $\G_{n, (p_2,\dots,p_K)}$}

For a multivariate normal approximation of the  joint vertex degrees $\bY = (Y_v, v \in \cM) \in \R^{M}$ we first obtain the covariance matrix.  

\begin{Lemma}\label{lem:cov}
Let $\bY \in \R^{M}$ be the vector of the joint vertex degrees in $\G_{n,(p_2,\dots,p_K)}$, i.e., $Y_v$ denotes the degree of vertex $v \in \cM$.
The covariance matrix $\Sigma = {\rm Cov} ( \bY)$ has entries 
\begin{align}
  \sigma_{\rm diag} :=  \Sigma_{v,v}&=  \sum_{k=2}^{{K}} {{n-1} \choose {k-1}} p_k (1-p_k)\\
   \sigma_{\rm off} :=   \Sigma_{u,v} &=    \sum_{k=2}^{{K}}  {{n-2} \choose {k-2}} p_k (1-p_k) 
   \qquad  \mbox{ for } u \ne v.  \nonumber 
\end{align}
\end{Lemma}
\noindent Its eigenvalues are $ {\sigma_{\rm diag} + (M-1) \sigma_{\rm off}} $, with multiplicity 1,  and 
${ \sigma_{\rm diag} - \sigma_{\rm off}} $, with multiplicity $M-1$.
If $p_k = p$ for all $k$ and if $M \ll n $ and $K=n$ then $\sigma_{\rm diag} = p(1-p) (2^{n-1} -1) $ and 
$\sigma_{\rm off} = p(1-p) 2^{n-2}$; the variance contributions are large.

For $M \ge 2$ we let
\begin{equation}\label{eq:tau}
\tau = \frac{\sigma_{\rm diag} - \sigma_{\rm off}}{\sigma_{\rm diag} + (M-1) \sigma_{\rm off}} \in (0,1]  
\end{equation}
denote the ratio of the smallest and the largest eigenvalue of $\Sigma$. For $M=1$ we set $\tau=1$. 
Then, for $K < n$, 
$$ \tau \ge \frac{n-K}{n-1 + (M-1)(K-1)}$$
and for $K \le n/2$ we obtain that $\tau \ge \frac{1}{2K}.$

For a multivariate normal approximation we employ \citep[Theorem 3.2]{temvcinas2024multivariate}, noting that we have a dissociated sum representation. The application requires a new smoothing lemma, which we derive; see the proof of Theorem \ref{th:normalp} in the appendix. 
We recall that $\E Y_v= \lambda $ given in \eqref{eq:mean} and note that $\Sigma$ in Lemma \ref{lem:cov} is positive semidefinite. 
In principle a natural object for a normal approximation is $\Sigma^{-1/2} (\bY - \lambda {\mathbf{1}})$, {when $\Sigma$ is positive definite,} but that transformation of $\bY$ obstructs the dependence structure. Instead we divide all $Y_v$ by
$$ \sigma= \sigma(n, M) := \sqrt{ \sigma_{\rm diag} + (M-1) \sigma_{\rm off}}.$$
\begin{Theorem}\label{th:normalp}
   Let $\bY \in \R^{M}$ be the vector of the joint vertex degrees in $\G_{n,(p_2,\dots,p_K)}$ and let $ W = \frac{1}{\sigma}( \bY - \lambda {\mathbf{1}})$ be the scaled vector of centred counts.  Let $\mathcal{K}$ be a class of convex sets in $\R^M$ and let  $Z \sim {\mathcal{MVN}}({0}, {\text{\rm Id}_{M \times M}})$ be a standard $M$-dimensional Gaussian vector. {Assume that $\Sigma$ is positive definite.}  Then \[
   d_{{\mathcal C}}( \LL (W), \LL (\sigma^{-1} \Sigma^{\frac{1}{2}}Z) )  
\le
\epsilon (n, M, (p_2, \ldots, p_K)) := \min\left\{1,
 \frac{2^{13/2}}{3^{3/4}}M^{9/16}
 B_{{\rm norm}}^{1/4}\sigma^{-3/4}\tau^{-3/8}\right\} ,
\]
where 
\begin{align*}
B_{{\rm norm}} &= \frac{4 M}{3} \sum_{k=2}^{K}  \left( {{n-1} \choose {k-1}} + 2 M {{n-2} \choose {k-2}} + M^2 {{n-3} \choose {k-3}} \right) p_k (1 - p_k)
\end{align*}
and  ${n-3 \choose k-3}:=0$ if $k=2$. 
\end{Theorem}

We note that, for the approximating normal distribution, $\sigma^{-2} \Sigma$ has largest eigenvalue 1 by construction. However it is theoretically possible that the covariance matrix $\Sigma$ becomes degenerate as $n \rightarrow \infty$.  We also note that a sharper bound could be obtained using Theorem 1.1 in \cite{raivc2019multivariate}, but at the cost of having to transform the degree vector in a way that the individual contributions cannot easily be read off without inverting the covariance matrix.

\begin{Remark}
     We 
have 
\begin{align*}
B_{{\rm norm}} & =\frac{4 M}{3}\left[\sigma^2(n, M)+(M+1) \sum_{k=2}^{K} \binom{n-2}{k-2} p_k\left(1-p_k\right)+M^2 \sum_{k=3}^{K}  \binom{n-3}{k-3} p_k\left(1-p_k\right)\right] \\
& 
\le (M+1)^2 \frac{4 M}{3} \sigma^2(n, M). 
\end{align*}
and hence 
$ {M^{9/4}} B_{{\rm norm}} \sigma^{-3} (n,M) 
\le \frac{4 (M+1)^{21/4}}{\sigma(n,M)}$.
So, if $\sigma(n,M)$ is  $\Theta(n)$, {if there is a $\tau_0>0$ such that $\tau = \tau(n) \ge \tau_0 $,} and if $M (n) = {o(n^{\frac{4}{21}})}$ then the bound in Theorem \ref{th:normalp} 
will be small for large $n$. 
\end{Remark}

\begin{Example}\label{ex:test}{Application: assessing the fit to a $\G_{n,(p_2,\dots,p_K)}$ model} 

As in Example \ref{ex:degreesequence}, if we think of $\cM$ as a sample of size $M$ then we may want to use the sample to assess the null hypothesis that the $\G_{n,(p_2,\dots,p_K)}$ model is appropriate for the data. For this test we centre and scale the observed degrees for $v \in \cM$ to obtain $\tX_v, v \in \cM$. The convex set distance is generally difficult to compute; here we use a scalar statistic, 
namely
$$T (X) = 
\max_{v \in \cM} \tX_v.$$
Using the maximum is motivated by the observation that in many real networks there are some vertices of very high degree, violating the homogeneity assumption which is implicit in   $\G_{n,(p_2,\dots,p_K)}$, giving the distribution $\PP_0$ under the null hypothesis. The maximum is real-valued, and in this case the convex set distance is equivalent to the Kolmogorov distance; indeed, $d_{KS} \le d_{\mathcal{C}} \le 2 d_{KS}$, see Section 2.1 in  \cite{meliot2022gaussian}. 
As the pre-image of a sublevel set under the maximum function is again a convex set, we can apply Theorem \ref{th:normalp} with 
${\tilde Z} = \sigma^{-1} \Sigma^{1/2} Z
$ to obtain that  
\begin{align*}
 \sup_{x \in \R } |   \PP_0 (T(\tX) \le x) - \PP (T({\tilde Z}) \le x)| \le
 {\epsilon (n, M, (p_2, \ldots, p_K))}
.
\end{align*}
While there is no closed form for the distribution of the maximum of correlated normals, we can gauge the null hypothesis by carrying out many simulations from the correlated normals and assess the distance to the observed maximum. 

Alternatively one could exploit the bounds in Theorem \ref{th:normalp} as an index for projection pursuit; see \cite{jee2009projection} for an overview.
\end{Example}

\section{Experimental Results}\label{sec:simulation}

\subsection{Simulation results}
To illustrate the theoretical results we simulate Erd\"{o}s-R\'{e}nyi   hypergraphs and compare the observed  distribution distance of a collection of vertex degrees with the theoretical Stein bounds in Theorems \ref{thm:jointp} and \ref{th:normalp}. 
We choose as maximal hyperedge order $K=3$ for these experiments, while varying the hypergraph size $n$ and the hyperedge probabilities $p_k, k=2, 3$. 
 
\subsubsection{Total variation distance to  multivariate Poisson
} 

\begin{figure}[!t]
    \centering
    \subfloat[
     Varying hypergraph size $n$.
    ]{\includegraphics[width=0.48\textwidth]{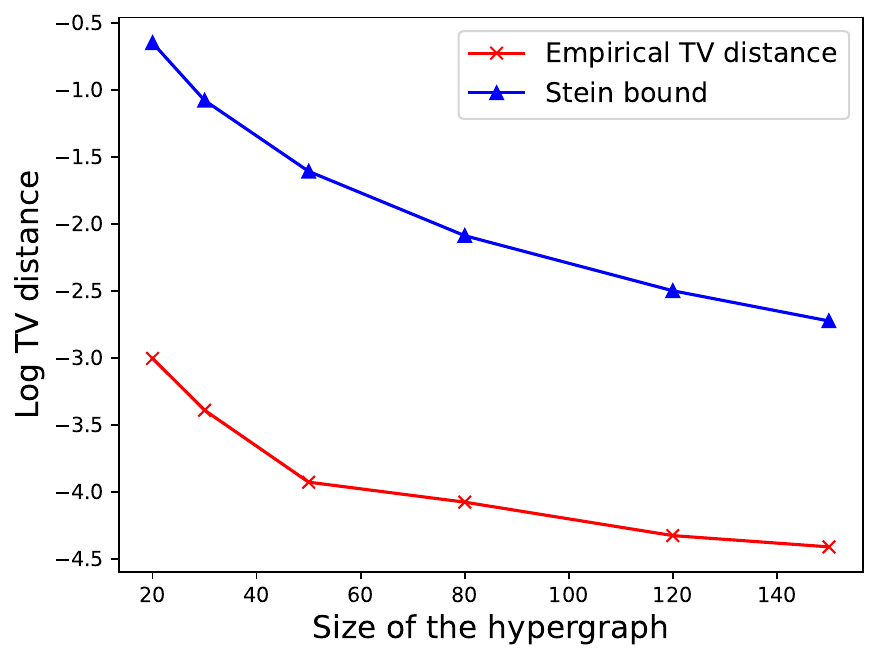}\label{fig:poisson_approx_p}}
    \subfloat[Varying $c$ in $p_k = \frac{c}{{n-1\choose k-1}}$, for 
    $n=100$.
    ]{\includegraphics[width=0.48\textwidth]{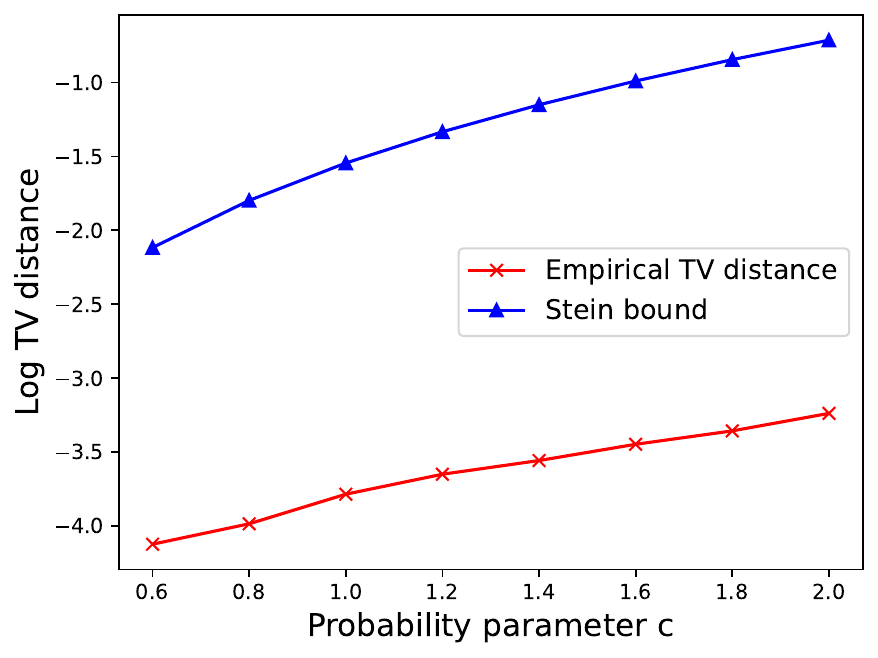}\label{fig:poisson_approx_n}}\par\medskip
 \subfloat[Degree sequence application as in Example \ref{ex:degreesequence} .
 ]{\includegraphics[width=0.48\textwidth]{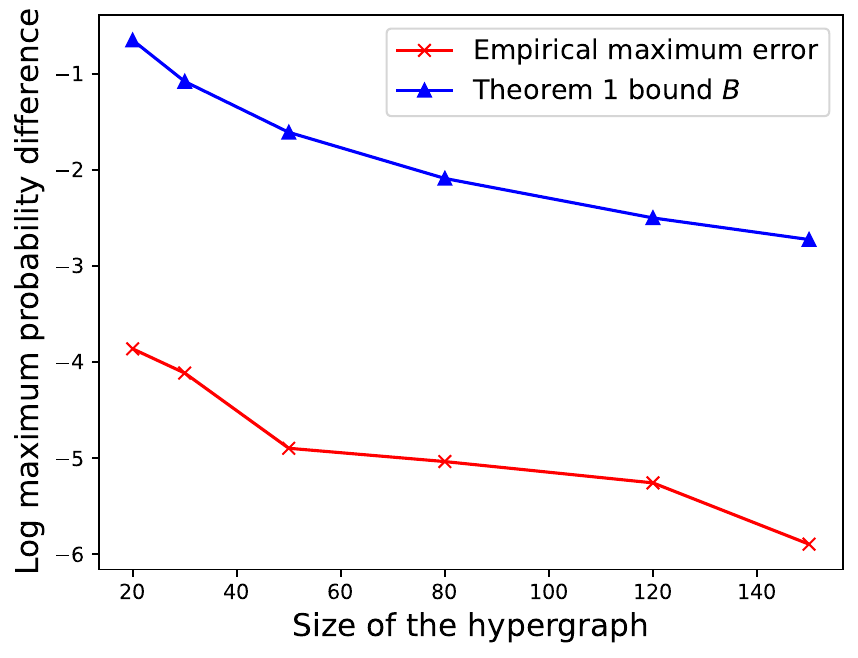}\label{fig:poisson_application}}
 \caption{Panels (a), (b): Bounds on the total variation distance to a multivariate Poisson distribution with $M=3, K=3, p_k = \frac{{c}}{{n-1\choose k-1}}$, and 100,000 trials. Panel (a): $c = 0.5$, hypergraph size $n$ varies; Panel (b): $n=100$, and $c$ varies; Panel (c): application to a degree sequence. Here $B= B_P$ is as in Theorem \ref{thm:jointp}.}
    \label{fig:poisson_approx}
\end{figure}

In line with Item 4 in 
Example \ref{ex:Poisson}  
we set $M=3,
p_k = \frac{0.5}{{n-1 \choose k-1}}$, for $k=2,3$ and use the log scale for the outcome. 
\Cref{fig:poisson_approx_p}   shows that the Stein bound and the observed total variation distance decrease to zero with increasing hypergraph size $n$. Moreover, the approximation quality decreases with increasing probability. We also note that the theoretical bound and the {observed}  bound behave similarly; it may be possible to improve on the constant in the theoretical bound.

\Cref{fig:poisson_application} compares $$\max_{0\leq \ell_v\leq M}|\PP ( Q_{k_v} \le \ell_v / M, v \in J) - \PP (d_{k_v} (\cM) \le \ell_v / M , v\in J )|$$ in Example \ref{ex:degreesequence}, with  degree pair $(k_1,k_2)=(0,1)$,  to  the bound 
$B$ in Theorem \ref{thm:jointp} for increasing hypergraph size $n$. 
The blue line is the theoretical upper bound $B$, and the red line is actual observed  probability difference in the simulation.  The probability difference follows a similar trend as the theoretical bound.

\subsubsection{Convex set distance to  multivariate normal
}

\begin{figure}[!tp]
    \centering
    \subfloat[
   $M=1$ and $p_k = 0.1$.]{\includegraphics[width=0.48\textwidth]{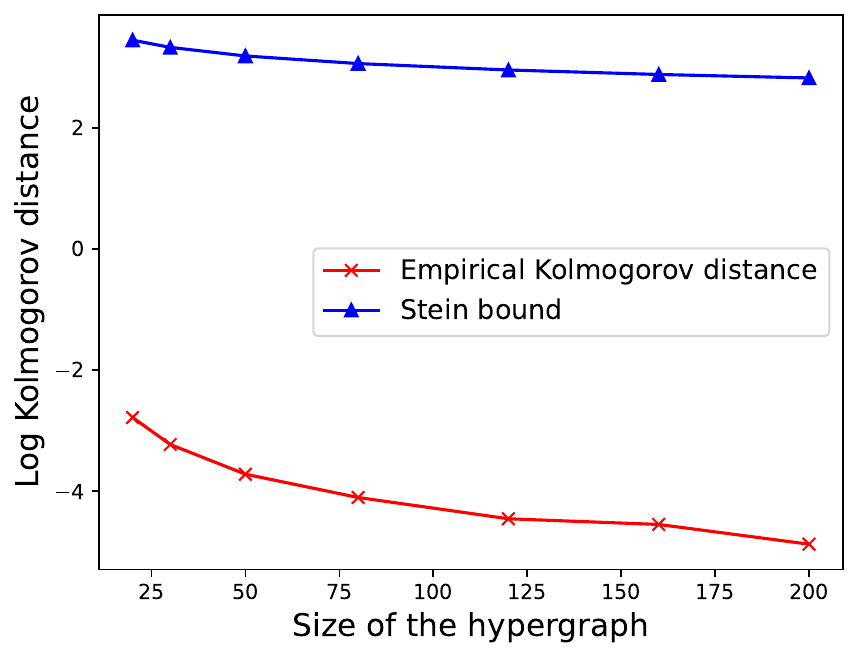}\label{fig:normal_approx_changen}}
    \subfloat[
    $M=1$,  and $n = 100$.]{\includegraphics[width=0.48\textwidth]{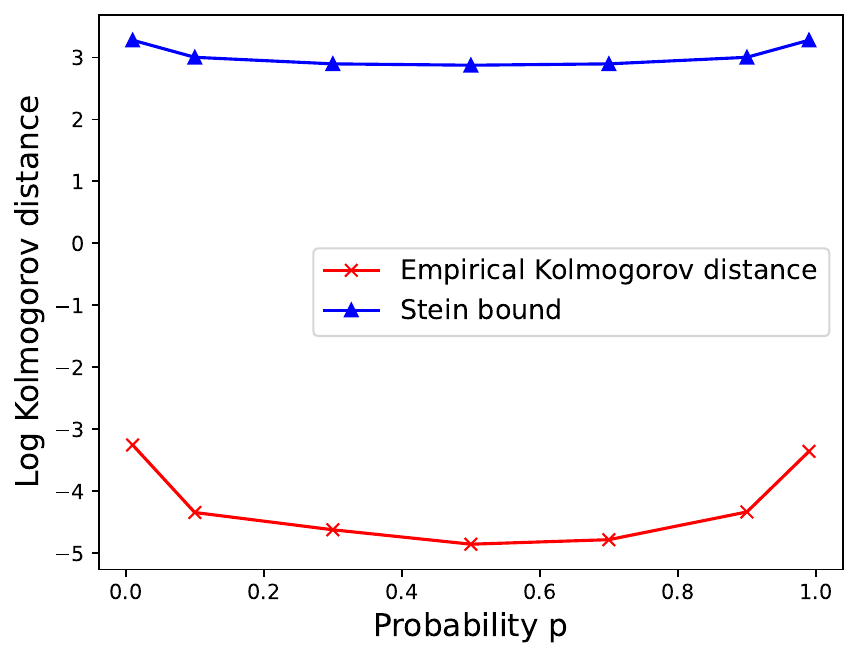}\label{fig:normal_approx_chengep}}\par\medskip
    \subfloat[
    $M=10$, and  $\delta=0.4$.]{\includegraphics[width=0.48\textwidth]{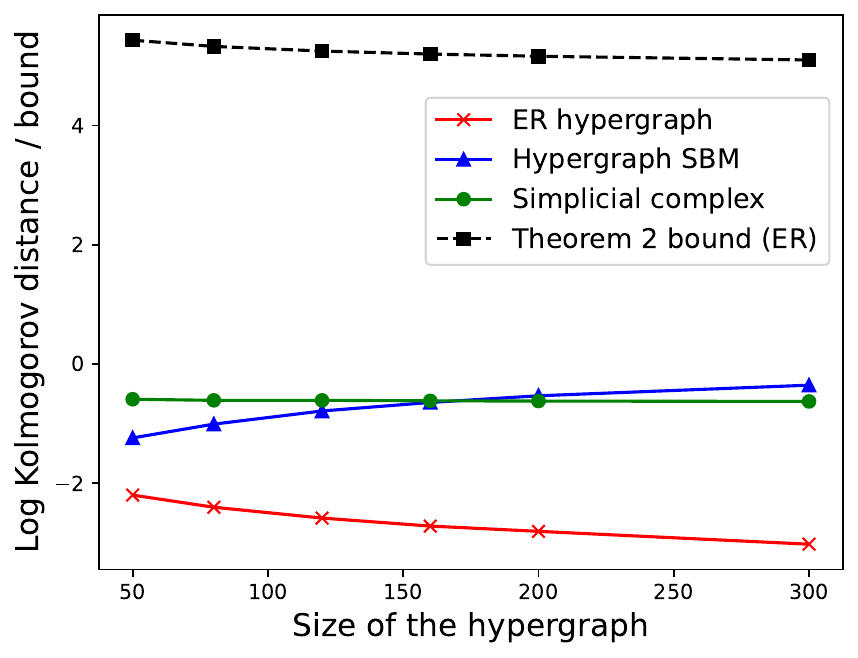}\label{fig:normal_test}}
    \subfloat[
    $M=10$, and $n=100$.]{\includegraphics[width=0.48\textwidth]{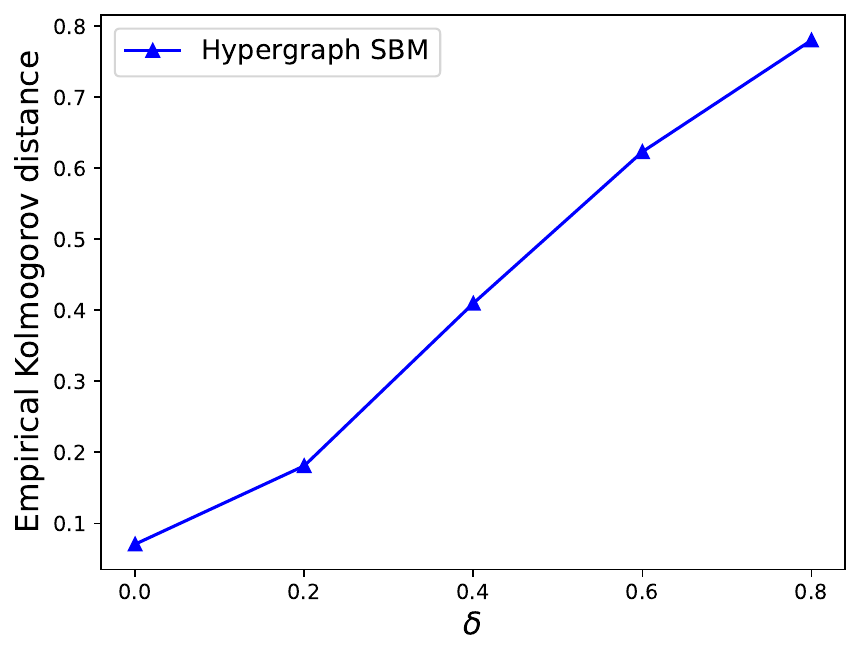}\label{fig:delta_test}}
    \par\medskip
    \subfloat[
    $M=10$, and  $\delta=0.4$.]{\includegraphics[width=0.48\textwidth]{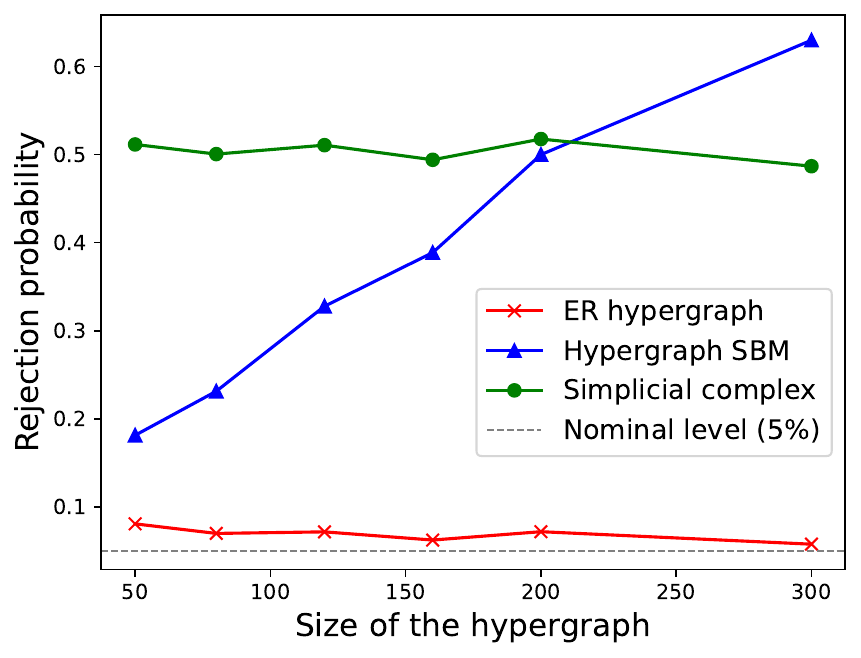}\label{fig:normal_rej}}
    \subfloat[
    $M=10$, and $n=100$.]{\includegraphics[width=0.48\textwidth]{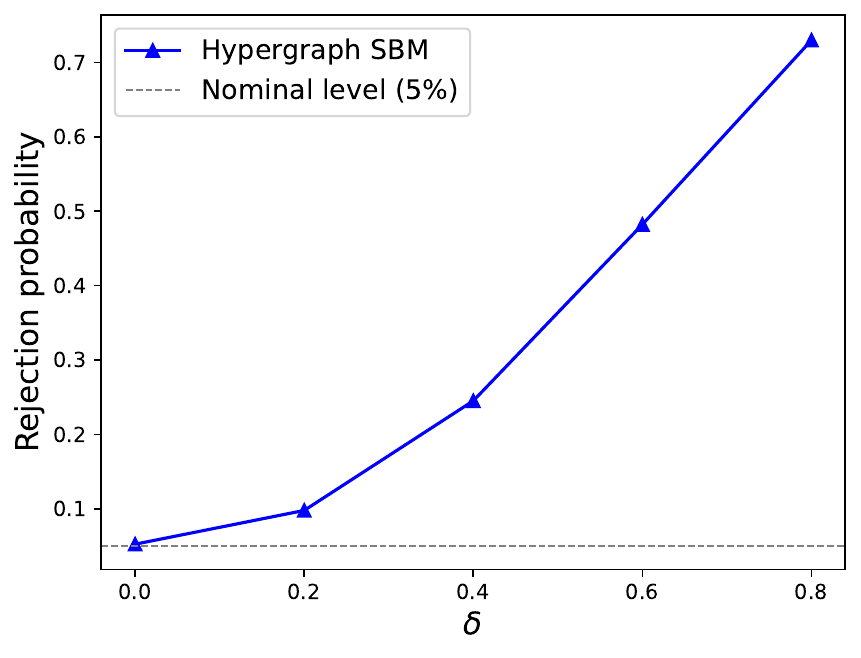}\label{fig:delta_rej}}
    \caption{Comparison to the multivariate normal distribution, showing the uncapped bound.
    For panels (a) and (b) we take $p_2 = p_3 = p_k$, for (c-f) we take $p_2 = 0.2, p_3 = 1 - (1-p_2)^{1/(n-2)};$  
  panels (a) and (b) illustrate  the  multivariate normal approximation for the degree vector, panels (c) and (d) show the Kolmogorov distance for the maximum degree, 
   whereas panels (e) and (f) 
   illustrate the application in Example \ref{ex:test} for different hypergraph models. For clarity we omit the truncation by 1 in the theoretical bounds here, and instead report the uncapped bounds. Panels (a) and (b) use 100,000 trials; panels (c-f) use 1,000,000 trials. All experiments use $K=3$. 
   }
    \label{fig:normal_approx}
\end{figure}

To illustrate the bounds on the convex set distance to multivariate normal 
in Theorem \ref{th:normalp}, we fix the hyperedge probabilities $p_2=p_3=p_k$ for $k=2,3$.
\Cref{fig:normal_approx} shows the theoretical Stein bound in Theorem~2. In the first two figures, we choose $M=1$ and show  the Kolmogorov distance, which is easier to calculate and equivalent to the  convex set distance.

In \Cref{fig:normal_approx_changen}, both the theoretical Stein bound and the observed  Kolmogorov distance decrease as $n$ increases for $p_2=p_3=0.1$. In \Cref{fig:normal_approx_chengep}, we fix $n=100$ and change the $p_k$. Both curves first decrease and then increase, reaching their minimum at $p_k=0.5$. 

Computing the supremum over all convex sets becomes difficult for $M>1$. We therefore illustrate the performance of the bound for multivariate case through Application~\ref{ex:test}.
We compare the Erd\"{o}s-R\'{e}nyi hypergraph model $\mathbb{G}_{n, (p_2,p_3)}$ with a hypergraph stochastic block model (hypergraph SBM)  (see {\cite[Definition 1.7]{dumitriu2026optimal}}) and a random simplicial complex (see {\cite[Definition 1.7]{cooley2020vanishing}})  under $
M=10, p_2=p=0.2, p_3=q_n=1-(1-p)^{1/(n-2)}.
$

In the hypergraph SBM we keep $K=3$ and take two 
equally sized communities $A$ and $B$;  
then
$
w_k=\frac{2\binom{n/2}{k}}{\binom nk}
$ 
is the proportion of potential $k$-hyperedges lying entirely within either community, $k=2,3$. 
The probabilities of a $k$-hyperedge within the first community, within the second community, and between communities are taken to be, respectively,
\[
p_k^A=p_k[1+\eta(1-w_k)+\delta],\qquad
p_k^B=p_k[1+\eta(1-w_k)-\delta],\qquad
p_k^{\mathrm{out}}=p_k(1-\eta w_k),
\]
where the parameter $\eta$ controls the preference for within-community hyperedges, and $\delta$ introduces a difference in expected degrees between the two communities.  
If we pick a potential hyperedge of size $k$ at random, the probability that it is included in the hypergraph SBM is $\frac{w_k}{2} p_k^A +\frac{w_k}{2} p_k^B + (1-w_k) p_k^{\mathrm{out}} = p_k$;
the setting  keeps the expected number of hyperedges of size $k$  
equal to that of $\mathbb{G}_{n, (p_2,p_3)}$.
 In our simulations we take  $\eta=0.8$ and $\delta=0.4$. 

For the simplicial complex, we independently include each size-$3$ hyperedges with probability $q_n$ and then deterministically include all of their subsets as hyperedges. By construction, for a given vertex pair, there are $n-2$ possible hyperedges of size $3$ that contain it, and the probability that this construction does not include it as $2$-hyperedge is $(1-q_n)^{n-2} = 1-p$. Hence, by construction, the probability $p_2$ of a $2$-hyperedge in this construction is $p_2 = p$.  Thus, every vertex has the same expected hyperedge degree as under the Erd\"{o}s-R\'{e}nyi hypergraph model, and there is no vertex heterogeneity.  However, this model has dependence between  
hyperedge occurrences, and an increasing hyperedge degree variance.

In \Cref{fig:normal_test}, as $n$ increases, both the Stein bound from Theorem \ref{th:normalp} and the Kolmogorov distance decrease under $\mathbb{G}_{n, (p_2,p_3)}$. However, the Kolmogorov distance increases under the hypergraph SBM due to degree heterogeneity. Under the simplicial complex, the distances remain  nearly unchanged when the hypergraph size changes, but they remain large compared to $\mathbb{G}_{n, (p_2,p_3)}$.  \Cref{fig:delta_test} shows that the Kolmogorov distance increases as  $\delta$ in the {hypergraph} SBM increases; here  $n=100$ is fixed.
To evaluate the performance of the test, we use the $95\%$ quantile of the maximum of multivariate normal distribution as the critical value, corresponding to a significance level of $5\%$. In \Cref{fig:normal_rej}, the $\mathbb{G}_{n, (p_2,p_3)}$ rejection rate is slightly above $\alpha$ but decreases towards it as $n$ grows. As $n$ increases, the rejection rate increases under the hypergraph SBM as $n$ increases,  under fixed $\delta=0.4$, and also  as $\delta$ increases, under fixed $n=100$, see \Cref{fig:delta_rej}, but the rejection rate remains stable under the simplicial complex when $n$ increases.

\subsection{Real Data Applications}

As real-world applications of our approximation results, we consider the MADStat coauthorship data from \cite{ji2022co} and the  primary-school contact data from \cite{stehle2011high}. As the underlying distribution is unknown, we use the empirical distribution of the degrees instead. We highlight  that Theorems \ref{thm:jointp} and \ref{th:normalp} do not control the error which arises from using an empirical distribution instead of the true underlying distribution.

\subsubsection{MADStat Coauthorship Data}

\begin{table}[tbp]
\centering
\small
\caption{MADStat hyperedge counts and estimated probabilities,
$\hat p_k=m_k/\binom{n}{k}$, with $n=47\,311$.}
\label{tab:madstat-probabilities}
\begin{tabular}{@{}rrl| @{\qquad}rrl@{}}
\toprule
$k$ & $m_k$ & $\hat p_k$
& $k$ & $m_k$ & $\hat p_k$ \\
\midrule
 2 & 24,295 & $2.171\times10^{-5}$
   & 18 & 1 & $4.554\times10^{-69}$ \\
 3 & 13,941 & $7.899\times10^{-10}$
   & 19 & 1 & $1.830\times10^{-72}$ \\
 4 & 4,422 & $2.119\times10^{-14}$
   & 20 & 0 & 0 \\
 5 & 1,308 & $6.623\times10^{-19}$
   & 21 & 0 & 0 \\
 6 & 509 & $3.269\times10^{-23}$
   & 22 & 2 & $3.197\times10^{-82}$ \\
 7 & 201 & $1.910\times10^{-27}$
   & 23 & 0 & 0 \\
 8 & 86 & $1.382\times10^{-31}$
   & 24 & 1 & $3.946\times10^{-89}$ \\
 9 & 41 & $1.254\times10^{-35}$
   & 25 & 0 & 0 \\
10 & 25 & $1.616\times10^{-39}$
   & 26 & 0 & 0 \\
11 & 5 & $7.517\times10^{-44}$
   & 27 & 0 & 0 \\
12 & 4 & $1.526\times10^{-47}$
   & 28 & 0 & 0 \\
13 & 2 & $2.097\times10^{-51}$
   & 29 & 0 & 0 \\
14 & 2 & $6.206\times10^{-55}$
   & 30 & 1 & $1.509\times10^{-108}$ \\
15 & 0 & 0
   & 31 & 0 & 0 \\
16 & 2 & $6.658\times10^{-62}$
   & 32 & 1 & $6.698\times10^{-115}$ \\
17 & 0 & 0
   & 33 & 2 & $9.350\times10^{-118}$ \\
\bottomrule
\end{tabular}
\end{table}

\begin{table}[tbp]
\centering
\small
\caption{MADStat approximation results and theoretical bounds.
}
\label{tab:madstat-errors}
\begin{tabular}{@{}llrr@{}}
\toprule
Data & Approximation Measure & 
Observed Distance & Theoretical bound \\
& & & assuming $\mathbb{G}_{n, (p_2,\cdots, {p_K})}$\\
\midrule
Full
& {$d_{TV}$}  for Multivariate Poisson ($M=3$)
& 0.78636 & 0.00075 \\
Full
& Kolmogorov for Normal ($M=1$)
& 0.52469 & 1\\
\midrule
$k\geq 6$
& {$d_{TV}$} for Multivariate Poisson ($M=3$)
& 0.08707 & 0.00011 \\
$k\geq 6$
& Probability Difference for Application \ref{ex:degreesequence} 
& 0.03635 & 0.00011 \\
$k\geq 6$
& Kolmogorov for Normal ($M=1$)
& 0.49172 & 1\\
\bottomrule
\end{tabular}
\end{table}

MADStat contains 83,331 articles by 47,311 authors published during 1975–2015 in 36 journals in statistics, probability, and machine learning. Each author set from a paper forms a hyperedge. We consider all hyperedges of size at least two and count each hyperedge only once if it occurs several times. Let ${\hat m}_k$ denote the number of hyperedges of size $k$ in the data set.
We estimate the 
probabilities for hyperedges of size $k$,  and the mean hyperedge degree of a vertex, 
$$\hat p_k=\frac{m_k}{\binom nk}, \quad \hat\mu=\sum_{k=2}^n{n-1\choose k-1}\hat{p}_k.$$
Table \ref{tab:madstat-probabilities} reports the hyperedge counts for different sizes and estimated probabilities. 

For a multivariate Poisson approximation, we randomly select 100 groups of three authors ($M=3$) and compare their {observed}  joint degree distribution with $\prod_{i=1}^3\operatorname{Poisson}(\hat\mu) $ using total variation distance. For multivariate normal approximation with $M=1$, we randomly select 100 authors ($M=1$), and we compute the Kolmogorov distance between the standardised hyperedge degrees of 100 randomly selected authors and $\sigma^{-1} \Sigma^{\frac{1}{2}}Z$.

For the full hypergraph, the total variation and Kolmogorov distances are 0.78636 and 0.52469, shown in Table \ref{tab:madstat-errors}. These two large distances indicate that both the multivariate Poisson and the  normal distribution do not fit the joint degree vectors. We note that the fitted hypergraph model $\mathbb{G}_{n, (p_2,\cdots, p_K)}$   is sparse and Theorem 1 gives a multivariate Poisson bound of 0.00075 under the fitted $\mathbb{G}_{n, (p_2,\cdots, p_K)}$ model. The observed degree variance, 24.25279, greatly exceeds the $\mathbb{G}_{n, (p_2,\cdots, p_K)}$ variance, 2.55395, at mean degree 2.55397. This difference suggests that the $\mathbb{G}_{n, (p_2,\cdots, p_K)}$ model hypergraph does not fit the observed hypergraph well.

We then study large collaborations by keeping only hyperedges with at least six authors. With this setting, the 
{observed} Kolmogorov distance is still large  at 0.49172, consistent with the low mean degree of 0.13084. On the other hand, the total variation distance decreases to 0.08707, showing a closer match to the multivariate Poisson distribution. 

Next we use Application \ref{ex:degreesequence} to compare the numbers of zero- and one-degree authors in each group to the multinomial distribution. The probability difference 
is 0.03635, showing a close match for these counts. We note that the theoretical multivariate Poisson bound is much smaller compared to the full model.  

To gauge how the sampling procedure may affect the observed distances, we selected 300 vertices randomly without replacement and split them into 100 disjoint sets of size 3. We use 
the fitted Erd\"os-R\'enyi hypergraph model with estimated probabilities and hyperedge size at least six and recorded the observed total variation distance. We repeated 
this experiment   100,000 times, generating a fresh hypergraph for each Monte Carlo replicate; the fitted parameters were held fixed across replicates.  
The central $95\%$ simulation interval is provided by its empirical 2.5th and 97.5th percentiles across the 100,000 Monte Carlo replicates. The observed  joint total variation  distance had mean 0.08231 and a $95\%$ interval [0.04184, 0.13867], compared with 0.08707 for the observed data. Thus, the total variation distance for the observed data is not unusual under the  fitted Erd\"os-R\'enyi hypergraph model.  
For Application \ref{ex:degreesequence}, the mean maximum joint cumulative probability difference was 0.04882, with a $95\%$ interval [0.01187, 0.11365], compared with the observed value 0.03635.  Again, the observed probability difference distance  is not unusual under the  fitted Erd\"os-R\'enyi hypergraph model.

\subsubsection{Primary School Contact Data}
Our second real-world data set is
a contact network dataset collected in a primary school in Lyon, France {\cite{stehle2011high}}. The data consist of recorded face-to-face interactions among 232 students and 10 teachers across 10 classes (one teacher per class), measured over two days in October 2009 using the SocioPatterns RFID infrastructure.

In our analysis, we focus on the higher-order interaction structure formed by these interactions. The original dataset includes only pairwise interactions. These interactions are recorded when two participants face each other within close distance (1 to 1.5 meters) during a short time window of 20 seconds. Therefore, if three interactions, A-B, B-C, and A-C, occur at the same timestamp, it is highly likely that students A, B, and C were engaged in a group interaction, given the tight spatial and temporal constraints. Therefore, to construct a hypergraph structure, if a set of student contacts forms a maximal clique of size $k$ at a given timestamp, we represent these contacts as size $k$ hyperedges. Smaller size hyperedges are not counted separately within the same time interval. The distribution of $k$-hyperedges and their estimated connection probabilities $\hat{p}_k$ for different values of $k$ is shown in Table 3. 

\begin{table}[tbp]
\centering
\small
\caption{Primary-school hyperedge counts and estimated connection
probabilities. }
\label{tab:school_edges}
\begin{tabular}{@{}rrr|rr@{}}
\toprule
& \multicolumn{2}{c}{Full school ($n=242$)}  
& \multicolumn{2}{c}{Class 5A ($n=22$)} \\
\cmidrule(lr){2-3}\cmidrule(lr){4-5}
$k$ & $m_k$ & $\hat p_k$ & $m_k$ & $\hat p_k$ \\
\midrule
2 & 7,748 & 0.26570 & 224 & 0.96970 \\
3 & 4,600 & 0.00197 & 219 & 0.14221 \\
4 & 347 & $2.489\times10^{-6}$ & 9 & 0.00123 \\
5 & 9 & $1.356\times10^{-9}$ & 0 & 0 \\
\bottomrule
\end{tabular}
\end{table}

\begin{table}[tbp] 
\centering
\small
\caption{Primary school approximation results and theoretical bounds. }
\label{tab:school_results}
\begin{tabular}{@{}llrr@{}}
\toprule
Data & Approximation Measure & 
{Observed} Distance & Theoretical bound  \\
\midrule
Full school & {$d_{TV}$} for Multivariate Poisson ($M=3$) & 0.99999 & 1\\
Full school & Kolmogorov for Normal ($M=1$) & 0.47609 & 1\\
Class 5A    & {$d_{TV}$}   for Multivariate Poisson ($M=3$)  & 0.99567 & 1\\
Class 5A    & Kolmogorov for Normal ($M=1$) & 0.19186 & 1\\
\bottomrule
\end{tabular}
\end{table}

We use 100 groups of three vertices $(M=3)$ for the multivariate Poisson approximation and 100 vertices $(M=1)$ for the multivariate normal approximation. The observed   total variation distance to multivariate Poisson, and Kolmogorov distance to normal, are 0.99999 and 0.47609, respectively; see Table \ref{tab:school_results}. Note that the hypergraph is relatively dense, with $\hat p_2=0.26570$. We also apply Application~\ref{ex:test} to 100 groups of three distinct vertices
($M=3$).  
We calculate the Monte Carlo $p$-values as $p_j=(1+\sum_{b=1}^{B}\mathbbm{1}(T_b^G\ge T_j))/(B+1)$ for the sampled statistic in group $j=1, \ldots, 100$,
where $T_j$ is its standardised maximum degree and $T_b^G$
are the maximum values of $B=1,000,000$ simulated correlated Gaussian vectors. We report the average $p$-value 
$\bar p=100^{-1}\sum_{j=1}^{100}p_j$,  and the rejection rate 
$R=100^{-1}\sum_{j=1}^{100}\mathbbm{1}(p_j\le0.05)$. 
We note that the average $p$-value is not itself an overall hypergraph
$p$-value. The 
observed Kolmogorov distance between maximum degree and normal distribution and rejection rate are large, 
shown in Table~\ref{tab:school_app2},
providing evidence against the fitted $\mathbb{G}_{n, (p_2,\cdots, {p_K})}$ model.  Similarly to what we saw in the  in the simulation results for the hypergraph SBM, these 
findings may partly reflect the community structure, which is not captured by the homogeneous $\mathbb{G}_{n, (p_2,\cdots, {p_K})}$ hypergraph model. We note here that the rejection rate is the proportion of rejected tests under 100 draws from the same hypergraph. 

\begin{table}[!t]
\centering
\small
\caption{Application~\ref{ex:test} for the primary-school data. Here ${\bar{p}}$ is the average Monte Carlo $p$-value and $R$ is the rejection rate.  
}
\label{tab:school_app2}
\begin{tabular}{@{}lrrc c @{}}
\toprule
Data & 
Observed Kolmogorov Distance & Theoretical bound & 
{${\bar{p}}$} 
& {$R$} 
\\
\midrule
Full school & $0.714776$ & $1$ &$0.15318$& $74\%$ \\
Class 5A    & $0.38669$ & $1$ &$0.27287$& $33\%$ \\
\bottomrule
\end{tabular}
\end{table}

We then consider all 22 students in Class 5A. For the multivariate normal approximation with $M=1$, we use each student’s degree exactly once. The joint total variation distance remains large at 0.99567, while the Kolmogorov distance decreases to 0.19186. 
The Application~\ref{ex:test} results in Table~\ref{tab:school_app2},
with a Kolmogorov distance $0.38669$ 
and rejection rate of $33\%$, also provide evidence against
the fitted $\mathbb{G}_{n, (p_2,\cdots, p_K)}$  model. The theoretical bounds 
equal one in both examples and thus provide no informative error control.

Again, to gauge how the sampling could affect the observed discrepancies, for Class 5A, we perform 100,000 Monte Carlo replicates under the fitted Erd\"os-R\'enyi hypergraph model, using the estimated hyperedge probabilities from the class 5A data. Each iteration uses all 22 vertex degrees for the multivariate normal comparison with $M=1$ and 100 randomly sampled groups of three distinct vertices for Application \ref{ex:test}. The observed Kolmogorov distances to the multivariate normal application have 
mean $0.2319$1, with empirical $95\%$ interval $[0.11019, 0.43392]$. For Application~3, the maximum-degree Kolmogorov distance has mean $0.32943$ and {empirical} $95\%$ interval $[0.15409, 0.60654]$. Thus, neither observed Kolmogorov distance  -- $0.19186$ for the degrees and $0.38669$ for Application \ref{ex:test} -- is   unusual under the fitted Erd\"os-R\'enyi hypergraph model. For the test in Application  \ref{ex:test}, we generating synthetic networks under the fitted  Erd\"os-R\'enyi hypergraph model and repeat the same test for the synthetic data. 
In our experiment
the Monte Carlo $p$-values 
have  
 average  $0.50131$ and a $95\%$ simulation interval $[0.22021, 0.79884]$; the observed average $p$-value $\bar{p} = 0.27287$ lies within this interval. 
However, in the simulated data the rejection rate of a test at level 5\% has mean $4.80\%$ and a  $95\%$ simulation 
interval $[0\%, 30\%]$, compared with $33\%$ for the observed data. 
The observed rejection rate exceeds the upper simulation limit, pointing against  
the Erd\"os-R\'enyi hypergraph model for Class~5A. This conclusion has to be qualified as our procedure is a plug-in Monte Carlo assessment; we only consider the empirical data distribution, and on a small sample.

\section{Discussion and Future Work}\label{sec:discussion}

This paper studies joint vertex-degree distributions in non-uniform Erd\"os-R\'enyi hypergraphs. Using Stein’s method, we establish explicit bounds for comparisons to multivariate Poisson under total variation distance and multivariate normal approximations under convex set distance. Two applications illustrate their use through a multinomial approximation for the joint degree counts among selected vertices, and a maximum-degree goodness-of-fit test under Kolmogorov distance. Simulation studies illustrate the theoretical results, while applications to coauthorship and contact data demonstrate their practical usefulness. The synthetic experiments illustrated the validity of the approximations, but also that there may be room for improvements in the constants appearing in the theoretical bounds. In the real-world examples, the main use was in  
assessing the $\mathbb{G}_{n, (p_2,\cdots, {p_K})}$ model for the data, to address whether 
there is more structure in the hypergraph than a $\mathbb{G}_{n, (p_2,\cdots, {p_K})}$  model accommodates. Both real-data examples indicate departures from the homogeneous $\mathbb{G}_{n, (p_2,\cdots, {p_K})}$ model; they also illustrate the importance of including the variation of the  empirical discrepancy statistics in the conclusions.

In future work, more general subgraph counts, including hypergraph triangles and other small motifs, could be studied in a similar way. 
Moreover, many objects of interest, such as motif counts, may occur in clumps; hence a compound Poisson approximation may be more appropriate; for an example in simple graphs see  for example \cite{picard2008assessing}. Hence, compound Poisson approximation bounds could be developed for hypergraph statistics. These approximating distributions for subgraph counts could potentially also support anomaly detection tasks in temporal hypergraphs.








\authorcontributions{
Conceptualization, {GR and WX}; methodology, {all authors}; software, {XZ and WX}; validation, {GR and XZ}; formal analysis, {GR and XZ}; investigation, {all authors};  data curation, {XZ}; writing---original draft preparation, {XZ and GR}; writing---review and editing, {all authors}; visualization, {XZ}; supervision, {GR and WX}; project administration, {GR and WX}; funding acquisition, {GR}. All authors have read and agreed to the published version of the manuscript.
}

\funding{{GR was funded in part by EPSRC grants EP/T018445/1, EP/X002195/1, EP/V056883/1, and 
EP/Y028872/1. 
WX was funded in part by EPSRC grant EP/T018445/1. 
}{XZ was funded in part by EPSRC DASS Programme grant  EP/Z531327/1 }
}


\institutionalreview{Not applicable.}



\dataavailability{ 
{The MADStat coauthorship data presented in the study are openly available in \cite{ji2016coauthorship} at \url{https://www.tracyke.net/MADStat.html}. The original 
primary-school contact data from  \cite{stehle2011high} used in the study are openly available  at \url{https://sociopatterns.org/datasets.html\#dataset-primary-school-cumulative-networks}. The processed primary-school contact data presented in the study are openly available in \cite{benson_simplicial_2018} at \url{https://www.cs.cornell.edu/~arb/data/contact-primary-school/} .}}

\acknowledgments{
{During the preparation of this manuscript/study, the authors used ChatGPT 6 for the purposes of proofreading and obtaining feedback{, as well as for help with the proof of Lemma \ref{lem:smoothing}}. The authors have reviewed and edited the output and take full responsibility for the content of this publication.}}

\conflictsofinterest{The authors declare no conflicts of interest.
The funders had no role in the design of the study; in the collection, analyses, or interpretation of data; in the writing of the manuscript; or in the decision to publish the~results.} 


\appendixtitles{no} 
\appendixstart
\appendix
\section{Technical Proofs}
\subsection{Proof of Theorem \ref{thm:jointp}}

\begin{proof} We note that 
$Y_v = \sum_{\alpha \in \Gamma_v} X_\alpha$, and that its mean $\lambda = \sum_{\alpha \in \Gamma_v} \mu_\alpha 
$ does not depend on $v$. 
   {First we note that by definition, the total variation distance is bounded by 1. Next,} we apply 
    {\citep[Theorem 10.A]{barbour1992poisson}} to  
     obtain that \begin{equation}
  d_{T V}\left(\mathcal{L}\left(\left\{Y_v\right\}_{v=1}^M\right), \prod_{v=1}^M \operatorname{Po}(\lambda) \right) \leq b_1 + b_2 
\end{equation}
where \begin{align*}
    b_1 &= {\sum_{\alpha \in \Gamma } (\mu_\alpha^2+\sum_{\beta \in B_\alpha \backslash\{\alpha\}} \mu_\alpha \mu_\beta )=} \sum_{\alpha \in \Gamma }\sum_{\beta \in B_\alpha} \mu_\alpha \mu_\beta ,\\
    b_2&={\sum_{\alpha\in\Gamma}\mathbb{E}\left[X_\alpha \sum_{\beta \in B_\alpha \backslash\{\alpha\}} X_\beta\right]}= \sum_{\alpha \in \Gamma } \sum_{\beta \in B_\alpha, \beta \ne \alpha } \E X_\alpha X_\beta . 
\end{align*}
The fourth term in \citep[Theorem 10.A]{barbour1992poisson} is zero due to the independence between distinct hyperedges. Next we calculate 
\begin{align*}
     b_1 &= \sum_{\alpha \in \Gamma } \sum_{\beta \in B_\alpha} \mu_\alpha \mu_\beta \\
     &= \sum_{v\in \cM } \sum_{u \in \cM} \sum_{e \in \cE: u \in e, v \in e} \left( \sum_{k=2}^n \mathbbm{1}(|e| =k) p_k \right)^2 \\ 
     &= \sum_{v\in \cM } \sum_{u\in \cM } \sum_{e \in \cE: u \in e, v \in e} \sum_{k=2}^n \mathbbm{1}(|e| =k) p_k^2 \\ 
     &\le M\sum_{k=2}^n {n-1 \choose k-1}p_k^2+M(M-1) \sum_{k=2}^n {{n-2} \choose {k-2}} p_k^2. 
\end{align*}

Here we used that each $e$ has exactly one size, and hence the sum in the parenthesis has exactly one non-zero element. 
Further, with \eqref{eq:dep}, 
\begin{align*}
b_2&= \sum_{\alpha \in \Gamma } \sum_{\beta \in B_\alpha, \beta \ne \alpha } \E X_\alpha X_\beta \\
&= \sum_{\alpha \in \Gamma } \sum_{\beta \in B_\alpha, \beta \ne \alpha } \mu_{\alpha}\\
     &= \sum_{v\in\cM}  \sum_{u\ne v}  
   \sum_{e \in \cE: u, v\in e} \sum_{k=2}^n  \mathbbm{1}( |e| =k) p_k \\
      &\le M(M-1) \sum_{k=2}^n  {{n-2}\choose {k-2}}  p_k . 
\end{align*}
{Using that $p_k=0$ for $k > K$ gives the assertion.}
\end{proof}

\subsection{Proof of Lemma \ref{lem:cov}}

\begin{proof} 
With the notation from Theorem \ref{thm:jointp}, we have $Y_v = \sum_{\alpha \in \Gamma_v} X_\alpha$ so that 
{\begin{align*}
 \Sigma_{v,v} &=   {\rm Var} (Y_v) = {\rm Var} (\sum_{\alpha \in \Gamma_v} X_\alpha) \\
 &= \sum_{\alpha {=(v, e)} \in \Gamma_v}  {\rm Var} (X_\alpha) + 
 \sum_{\alpha{=(v, e)} \in \Gamma_v} \sum_{\beta {=(v, f)} \in \Gamma_v, {f \ne e} } {\rm Cov} (X_\alpha, X_\beta). 
\end{align*}
As the hyperedge indicators for $e$ and $f \ne e$ are independent, the second term vanishes. Moreover, 
for $\alpha = (v, e) \in \Gamma_v$
\begin{align} \label{eq:sigmaalpha}
 \sigma_\alpha :=   {\rm Var} X_\alpha &=
    \sum_{k=2}^n {\mathbbm{1} } ( |e| =k) p_k (1-p_k) 
\end{align} 
and hence 
\begin{align*}
 \Sigma_{v,v} &= \sum_{\alpha \in \Gamma_v}   {\rm Var} X_\alpha =
    \sum_{k=2}^n {{n-1} \choose {k-1}} p_k (1-p_k) . 
\end{align*} 
Similarly, for $\beta = (u,e) \in \Gamma_u$ with $u \ne v,$
} 
{\begin{align*}
\Sigma_{u, v} & =\operatorname{Cov}\left(Y_v, Y_u\right) 
 =\sum_{\alpha =(v,e) \in \Gamma_v} \sum_{\beta =(u,f) \in \Gamma_u} \operatorname{Cov}\left(X_\alpha, X_\beta\right) \\
& =\sum_{\alpha =(v,e) \in \Gamma_v} \sum_{\beta \in B_\alpha \cap \Gamma_u} \operatorname{Cov}\left(X_\alpha, X_\beta\right)\\
& =\sum_{\alpha =(v,e) \in \Gamma_v} \sum_{\beta =(u,e) \in B_\alpha \cap \Gamma_u} \operatorname{Cov}\left(X_\alpha, X_\beta\right)\\
&=\sum_{k=2}^n {\sum_{e \in \cE}} {\mathbbm{1} } ( |e| =k; u, v \in e) p_k (1-p_k)\\
&=\sum_{k=2}^{n} {{n-2} \choose {k-2}} p_k (1-p_k)
.
\end{align*}
}
{Using that $p_k=0$ for $k > K$ gives the assertion.} 
\end{proof}

\subsection{Proof of Theorem \ref{th:normalp}}

{In order to prove Theorem \ref{th:normalp} we first have to amend the smoothing argument leading Theorem 3.3 in \cite{temvcinas2024multivariate}. Let \(U\in\R^d\) be centered, with positive definite covariance matrix
\(\Sigma\). Set  
\begin{equation}\label{eq:normalization}
 \sigma^2=\lambda_{\max}(\Sigma),\qquad
 \tau=\frac{\lambda_{\min}(\Sigma)}{\lambda_{\max}(\Sigma)},\qquad
 W=\frac{U}{\sigma},\qquad V=\frac{\Sigma}{\sigma^2}.
\end{equation}
Thus $0<\tau\le  1$  and $\tau I_d\preceq V\preceq I_d$, where $I_d$ is the $d$-dimensional identity matrix. Here \(A\preceq B\) means that \(B-A\) is positive semidefinite. 
Let \(Z\sim N_d(0, I_d)\) and
\[
 G=\sigma^{-1}\Sigma^{1/2}Z\sim N_d(0,V).
\]  We say $h \in \testzero$ iff $h$ is three times partially differentiable with third partial derivatives being Lipschitz and bounded.
For a function $h \in \testzero$  we write $
 |h|_3=\max_{1\leq i,j,k\leq d}\|\partial_{ijk}h\|_\infty.
$  Then we show the following smoothing lemma.
}

\begin{Lemma}\label{lem:smoothing}
 {Let $ W=\frac{U}{\sigma}$ be as in \eqref{eq:normalization}. Assume that a finite \(\beta\geq0\) satisfies
\begin{equation}\label{eq:smooth-input}
 |\E h(W)-\E h(G)|\leq\beta\,|h|_3
\end{equation}
for every $h \in \testzero$. Then 
\begin{equation}\label{eq:main-bound}
d_{{\mathcal C}} (\LL (W),\LL (G))
 \le \min\left\{1,
 \frac{2^{13/2}}{3^{3/4}}\,
d^{9/16}\beta^{1/4}\tau^{-3/8}\right\}. 
\end{equation}
}
\end{Lemma}

\begin{proof} 
  {Let ${\mathcal C}_d$ denote the class of convex Borel subsets of $\R^d$, and let $\varepsilon>0$. For $A \in {\mathcal C}_d$ let $ {\rm dist} (x,A) = \inf_{y \in A} || x - y||_2$
  and set 
\[
 A^\varepsilon=\{x:{\rm dist} (x,A)\le \varepsilon\},\qquad
 A^{-\varepsilon}=\{x\in A:\overline{B}(x,\varepsilon)\subseteq A\},
\]
where \(\overline{B}(x,\varepsilon)\) is the closed Euclidean ball of radius $\epsilon$ around $x$. 
For $Z_0\sim N_d(0,I_d)$, we have 
\begin{equation}\label{eq:standard-boundary}
 \sup_{A\in {\mathcal C}_d }\max\bigl\{
 \PP (Z_0\in A^\varepsilon\setminus A),
 \PP (Z_0\in A\setminus A^{-\varepsilon})\bigr\}
 \leq 4d^{1/4}\varepsilon; 
\end{equation}
see Proposition~1.1 and the subsequent discussion in \cite{raivc2019multivariate}.}

{Since $V\succeq\tau I_d$, we have the representation 
\[
 G\overset{\mathrm d}=\sqrt{\tau}\,Z_0+R,
 \qquad R\sim N_d(0,V-\tau I_d),
\]
with $Z_0$ and $R$ independent; the covariance of $R$ may be degenerate.   Conditional on $R=r$, translate by $-r$  and divide by $\sqrt{\tau}$ . Convexity is preserved, and a boundary strip of width
$\varepsilon$ has width $\varepsilon/\sqrt{\tau}$ after this change
of coordinates. Applying \eqref{eq:standard-boundary} and  taking the expectation over 
$R$  gives
\begin{equation}\label{eq:boundary-tau}
 \sup_{A\in {\mathcal C}_d }\max\bigl\{
 \PP (G\in A^\varepsilon\setminus A),
 \PP (G\in A\setminus A^{-\varepsilon})\bigr\}
 \le 4d^{1/4}\tau^{-1/2}\varepsilon.
\end{equation}
Moreover the use of closed rather than open boundary strips does not change these
Gaussian probabilities, since $V$ is positive definite.
  } 

  {In the next step we approximate the indicator function of a convex set by a function $h$ which satisfies the assumptions of the lemma, as follows. For a bounded measurable function $f:\R^d\to\R$, let
\[
 (S_\delta f)(x)
 =\frac{1}{(2\delta)^d}\int_{[-\delta,\delta]^d} f(x+z)\,dz.
\]
Choose
\begin{equation*}
 \delta=\frac{\varepsilon}{16\sqrt d},\qquad
 h_{\varepsilon,A}=S_\delta^4 \, {\mathbbm 1}_{A^{\varepsilon/4}}.
\end{equation*}
Then \begin{equation}\label{eq:sandwich}
 0\leq h_{\varepsilon,A}\le 1,\qquad
 h_{\varepsilon,A}=1\text{ on }A,\qquad
 h_{\varepsilon,A}=0\text{ outside }A^\varepsilon.
\end{equation}
 Lemma 9 in \cite{gan2017dirichlet} gives that $h \in \testzero$ and 
\begin{equation}\label{eq:derivative-bound}
 |h_{\varepsilon,A}|_3
 \le \delta^{-3}
 =16^3d^{3/2}\varepsilon^{-3}.
\end{equation}
Moreover, Lemma 8 in \cite{gan2017dirichlet}, applied to \eqref{eq:sandwich}, states that
\begin{align}
d_{{\mathcal C}} (\LL (W),\LL (G)) 
 & \le
 \sup_{A\in {\mathcal C}_d}
 |\E h_{\varepsilon,A}(W)-\E h_{\varepsilon,A}(G)|
 \notag\\
 &\quad+
 \sup_{A\in {\mathcal C}_d}\max\bigl\{
 \PP (G\in A^\varepsilon\setminus A),
 \PP (G\in A\setminus A^{-\varepsilon})\bigr\}.
 \label{eq:smoothing-inequality}
\end{align}
Combining \eqref{eq:smooth-input}, \eqref{eq:boundary-tau} and
\eqref{eq:derivative-bound}, for every \(\varepsilon>0\) we obtain
\begin{equation}\label{eq:pre-optimization}
d_{{\mathcal C}} (\LL (W),\LL (G)) 
 \le 4d^{1/4}\tau^{-1/2}\varepsilon
       +16^3d^{3/2}\beta\varepsilon^{-3}. 
\end{equation}
  }
  {The last step is to optimise the smoothing parameter $\varepsilon$. For positive $a,b$,  differentiation yields
\[
 \inf_{\varepsilon>0}(a\varepsilon+b\varepsilon^{-3})
 =\frac{4}{3^{3/4}}a^{3/4}b^{1/4},
\]
and the infimum is achieved at $
 \varepsilon_*=(3b/a)^{1/4}.$ Set $a=4d^{1/4}\tau^{-1/2}$ and
$b=16^3d^{3/2}\beta$. When $\beta>0$, the resulting estimate is
\begin{align*}
 d_{{\mathcal C}} (\LL (W),\LL (G)) 
 &\le \frac{4}{3^{3/4}}
 (4d^{1/4}\tau^{-1/2})^{3/4}
 (16^3d^{3/2}\beta)^{1/4}\\
 &=\frac{2^{13/2}}{3^{3/4}}
 d^{9/16}\beta^{1/4}\tau^{-3/8}.
\end{align*} For $\beta=0$,  we let $\varepsilon\downarrow 0$ in
\eqref{eq:pre-optimization}. This gives the assertion. }
\end{proof} 

{Note that if $\beta = B/\sigma^3$ then \eqref{eq:main-bound} yields 
\begin{equation}\label{eq:raw-bound}
d_{{\mathcal C}} (\LL (W),\LL (G)) 
 \le \min\left\{1,
 \frac{2^{13/2}}{3^{3/4}}\,
 d^{9/16}B^{1/4}\sigma^{-3/4}\tau^{-3/8}\right\}.
\end{equation}} 
Now we can prove  Theorem \ref{th:normalp}. 
\begin{proof}
For $\alpha = (v,e)$ we set ${\tilde{X}}_\alpha = \frac{1}{\sigma}(X_\alpha - p_e)$, so that $\tX_\alpha$ has mean 0. 
    We apply \citep[Theorem 3.2]{temvcinas2024multivariate} ; matching the notation in that paper: $W_i = {\frac{1}{\sigma}} ( Y_v - \lambda), d = {M}, i = v, {\mathbbm{I}}_i = \Gamma_v, {\mathbbm{D}}_j(s) = B_\alpha {\cap\Gamma_j}, X_s = {\tilde{X}}_\alpha = {\frac{1}{\sigma}}(X_{(v,e)} - p_e)$.

    {Let $h \in \testzero.$ 
Then  \[ | \EE h(W) - \EE h( \sigma^{-1} \Sigma^{\frac{1}{2}}Z)| 
  \le  (B_1 + B_{2}) | h|_3,
\]} 
where 
\begin{align*}
B_{1} &=
\frac{1}{3}\sum_{\alpha \in \Gamma} \sum_{t, u \in B_\alpha}  \left(\frac{1}{2}\EE{
| \tX_\alpha \tX_t \tX_u| }+\EE {|\tX_\alpha \tX_t|} \EE {|\tX_u|}\right) \\
B_{2} &=  \frac{1}{3} \sum_{\alpha =(v,e) \in \Gamma} \sum_{t=(u,e) \in B_\alpha} \sum_{r=(w,e) \in B_\alpha \setminus B_t}  \left(\EE {|\tX_\alpha \tX_t \tX_v|} + \EE {|\tX_\alpha \tX_t|} \EE {|\tX_v|}\right). 
\end{align*}
Here we note that in this instance, $B_\alpha = B_t$ as the edge $e$ is fixed. Hence {$B_2=0$ and} we only need to bound  $B_1$. We have 
\begin{align*}
B_{1} &= B_{1,1} + B_{1,2}
\end{align*}
with 
\begin{align*} 
B_{1,1} &= 
\frac{1}{3} \sum_{v \in \cM} \sum_{e \in E: v \in e} \sum_{k=2}^{n} {\mathbbm{1}}  (|e| =k) \sum_{t, u \in \cM; t, u \in e}
\frac{1}{2}\EE{
| \tX_{(v,e)}  \tX_{(t,e)} \tX_{(u,e)}| }\\
B_{1,2} &= \frac{1}{3} \sum_{v \in \cM} \sum_{e \in E: v \in e} \sum_{k=2}^{n} {\mathbbm{1}}  (|e| =k) \sum_{t, u \in \cM; t, u \in e} \EE {|\tX_{(v,e)} \tX_{(t,e)}|} \EE {|\tX_{(u,e)}|}. 
\end{align*}
For $e$ with $|e| =k$ and $v \in e$, 
\begin{align*}
  \EE |\tX_{(v,e)}| &= \frac{1}{\sigma} 2 p_k (1-p_k).  
\end{align*}
Moreover, for $t \in B_\alpha$, the same edge $e$ as in $\alpha$ is used; hence 
\begin{align*}
\mathbb{E}\left|\tilde{X}_{(u, e)} \tilde{X}_{(t, e)}\right| 
& =\frac{1}{\sigma^2}\left[p_k\left(1-p_k\right)^2+\left(1-p_k\right) p_k^2\right] \\
& =\frac{p_k\left(1-p_k\right)}{\sigma^2}.
\end{align*}
Similarly, 
\begin{align*}
\mathbb{E}\left|\tilde{X}_{(v, e)} \tX_{(t, e)} \tX_{(u, e)}\right| 
& =\frac{p_k\left(1-p_k\right)\left(1-2 p_k+2 p_k^2\right)}{\sigma^3}
\end{align*}
Collecting these terms we have
\begin{align*} 
B_{1,1} 
&=  
\frac{1}{6\sigma^3} \sum_{v \in \cM} \sum_{e \in E: v \in e} \sum_{k=2}^{n} {\mathbbm{1}}  (|e| =k) \sum_{t, u \in \cM; t, u \in e}
p_k\left(1-p_k\right)\left(1-2 p_k+2 p_k^2\right)\\
&\le \frac{M}{6 \sigma^3} \sum_{k=2}^n\left[\binom{n-1}{k-1}+3(M-1)\binom{n-2}{k-2}+(M-1)(M-2)\binom{n-3}{k-3}\right] p_k\left(1-p_k\right) 
,\\
B_{1,2} &\le \frac{2 M}{3 \sigma^3} \sum_{k=2}^n\left[\binom{n-1}{k-1}+3(M-1)\binom{n-2}{k-2}+(M-1)(M-2)\binom{n-3}{k-3}\right]p_k^2\left(1-p_k\right)^2. 
\end{align*}
Using that $p_k(1-p_k) \le \frac14$
we obtain
\begin{align*}
    B_1&\le\frac{M}{3 \sigma^3} \sum_{k=2}^n\left[\binom{n-1}{k-1}+3(M-1)\binom{n-2}{k-2}+(M-1)(M-2)\binom{n-3}{k-3}\right]p_k\left(1-p_k\right) \\
    &{\le \frac{B_{{\rm norm}}}{\sigma^3}}
\end{align*}
{with $B_{{\rm norm}}$ as in the statement of Theorem \ref{th:normalp}.} 
{Using that $p_k=0$ for $k > K$ {and applying \eqref{eq:raw-bound}, { noting that the convex set distance is always bounded by $1$,} gives}  the assertion.}
\end{proof}
\reftitle{References}
\bibliography{main}

@article{cooley2020phase,
  title={Phase transition in cohomology groups of non-uniform random simplicial complexes},
  author={Cooley, Oliver and Del Giudice, Nicola and Kang, Mihyun and Spr{\"u}ssel, Philipp},
  journal={arXiv preprint arXiv:2005.07103},
  year={2020}
}

@article{picard2008assessing,
  title={Assessing the exceptionality of network motifs},
  author={Picard, Franck and Daudin, J-J and Koskas, Michel and Schbath, Sophie and Robin, Stephane},
  journal={Journal of Computational Biology},
  volume={15},
  number={1},
  pages={1--20},
  year={2008},
  publisher={SAGE Publications, Sage CA: Los Angeles, CA}
}

@article{jee2009projection,
  title={Projection pursuit},
  author={Jee, J Rodney},
  journal={Wiley Interdisciplinary Reviews: Computational Statistics},
  volume={1},
  number={2},
  pages={208--215},
  year={2009},
  publisher={Wiley Online Library}
}

@article{dumitriu2026optimal,
  title={Optimal and exact recovery on the general nonuniform hypergraph stochastic block model},
  author={Dumitriu, Ioana and Wang, Hai-Xiao},
  journal={The Annals of Statistics},
  volume={54},
  number={1},
  pages={48--73},
  year={2026},
  publisher={Institute of Mathematical Statistics}
}

@article{cooley2020vanishing,
  title={Vanishing of cohomology groups of random simplicial complexes},
  author={Cooley, Oliver and Giudice, Nicola Del and Kang, Mihyun and Spr{\"u}ssel, Philipp},
  journal={Random Structures \& Algorithms},
  volume={56},
  number={2},
  pages={461--500},
  year={2020},
  publisher={Wiley Online Library}
}

@article{stehle2011high,
  title={High-resolution measurements of face-to-face contact patterns in a primary school},
  author={Stehl{\'e}, Juliette and Voirin, Nicolas and Barrat, Alain and Cattuto, Ciro and Isella, Lorenzo and Pinton, Jean-Fran{\c{c}}ois and Quaggiotto, Marco and Van den Broeck, Wouter and R{\'e}gis, Corinne and Lina, Bruno and others},
  journal={PloS one},
  volume={6},
  number={8},
  pages={e23176},
  year={2011},
  publisher={Public Library of Science, San Francisco, USA}
}

@article{ji2022co,
  title={Co-citation and co-authorship networks of statisticians},
  author={Ji, Pengsheng and Jin, Jiashun and Ke, Zheng Tracy and Li, Wanshan},
  journal={Journal of Business \& Economic Statistics},
  volume={40},
  number={2},
  pages={469--485},
  year={2022},
  publisher={Taylor \& Francis}
}

@article{schawe2022higher,
  title={Higher order interactions destroy phase transitions in deffuant opinion dynamics model},
  author={Schawe, Hendrik and Hern{\'a}ndez, Laura},
  journal={Communications Physics},
  volume={5},
  number={1},
  pages={32},
  year={2022},
  publisher={Nature Publishing Group UK London}
}

@article{agostinelli2026higher,
  title={Higher-order dissimilarity measures for hypergraph comparison},
  author={Agostinelli, Cosimo and Mancastroppa, Marco and Barrat, Alain},
  journal={Journal of Complex Networks},
  volume={14},
  number={1},
  pages={cnaf048},
  year={2026},
  publisher={Oxford University Press}
}

@article{meliot2022gaussian,
  title={Gaussian approximations for random vectors},
  author={M{\'e}liot, Pierre-Lo{\"\i}c and Nikeghbali, Ashkan},
  journal={arXiv preprint arXiv:2203.01712},
  year={2022}
}

@article{temvcinas2024multivariate,
  title={Multivariate central limit theorems for random clique complexes},
  author={Tem{\v{c}}inas, Tadas and Nanda, Vidit and Reinert, Gesine},
  journal={Journal of Applied and Computational Topology},
  volume={8},
  number={6},
  pages={1837--1880},
  year={2024},
  publisher={Springer}
}

@book{barbour2026networks,
  title={Networks: {P}robability and {S}tatistics},
  author={Barbour, Andrew D and Reinert, Gesine},
  year={2026},
  publisher={Cambridge University Press}
}

@article{klimm2021hypergraphs,
  title={Hypergraphs for predicting essential genes using multiprotein complex data},
  author={Klimm, Florian and Deane, Charlotte M and Reinert, Gesine},
  journal={Journal of Complex Networks},
  volume={9},
  number={2},
  year={2021},
  publisher={Oxford University Press}
}

@book{chen2011normal,
  title={Normal {A}pproximation by {S}tein's {M}ethod},
  author={Chen, Louis HY and Goldstein, Larry and Shao, Qi-Man},
  volume={2},
  year={2011},
  publisher={Springer}
}

@book{barbour1992poisson,
  title={Poisson {A}pproximation},
  author={Barbour, Andrew D and Holst, Lars and Janson, Svante},
  year={1992},
  publisher={The Clarendon Press Oxford University Press}
}

@article{benson_simplicial_2018,
	title = {Simplicial closure and higher-order link prediction},
	volume = {115},
	url = {https://www.pnas.org/doi/abs/10.1073/pnas.1800683115},
	doi = {10.1073/pnas.1800683115},
	number = {48},
	urldate = {2025-03-22},
	journal = {Proceedings of the National Academy of Sciences},
	author = {Benson, Austin R. and Abebe, Rediet and Schaub, Michael T. and Jadbabaie, Ali and Kleinberg, Jon},
	month = {November},
	year = {2018},
	note = {Publisher: Proceedings of the National Academy of Sciences},
	pages = {E11221--E11230}
}

@article{gan2017dirichlet,
  title={Dirichlet approximation of equilibrium distributions in Cannings models with mutation},
  author={Gan, Han L and R{\"o}llin, Adrian and Ross, Nathan},
  journal={Advances in Applied Probability},
  volume={49},
  number={3},
  pages={927--959},
  year={2017},
  publisher={Cambridge University Press}
}

@article{raivc2019multivariate,
  title={A multivariate Berry--Esseen theorem with explicit constants},
  author={Rai{\v{c}}, Martin},
  journal={Bernoulli},
  volume={25},
  number={4A},
  pages={2824--2853},
  year={2019},
  publisher={JSTOR}
}

@inproceedings{stein1972bound,
  title={A bound for the error in the normal approximation to the distribution of a sum of dependent random variables},
  author={Stein, Charles},
  booktitle={Sixth Berkeley Symposium on Mathematical Statistics and Probability},
  volume={2},
  pages={583--603},
  year={1972},
  organization={University of California Press}
}

@article{darling2005structure,
  title={Structure of large random hypergraphs},
  author={Darling, Richard WR and Norris, James R},
  journal={The Annals of Applied Probability},
  volume={15},
  number={1A},
  pages={125--152},
  year={2005},
  publisher={Institute of Mathematical Statistics}
}

@article{chodrow2020configuration,
  title={Configuration models of random hypergraphs},
  author={Chodrow, Philip S},
  journal={Journal of Complex Networks},
  volume={8},
  number={3},
  year={2020},
  publisher={Oxford University Press}
}

@article{ji2016coauthorship,
  title={Coauthorship and citation networks for statisticians},
  author={Ji, Pengsheng and Jin, Jiashun},
  journal={The Annals of Applied Statistics},
  volume={10},
  number={4},
  pages={1779--1812},
  year={2016},
  publisher={Institute of Mathematical Statistics}
}

@article{austin2008exchangeable,
  title={On exchangeable random variables and the statistics of large graphs and hypergraphs},
  author={Austin, Tim},
  journal={Probability Surveys},
  volume={5},
  pages={80--145},
  year={2008},
  publisher={Institute of Mathematical Statistics and Bernoulli Society}
}

@article{ghoshal2009random,
  title={Random hypergraphs and their applications},
  author={Ghoshal, Gourab and Zlati{\'c}, Vinko and Caldarelli, Guido and Newman, Mark EJ},
  journal={Physical Review E},
  volume={79},
  number={6},
  pages={066118},
  year={2009},
  publisher={APS}
}

@article{zhou2006learning,
  title={Learning with hypergraphs: Clustering, classification, and embedding},
  author={Zhou, Dengyong and Huang, Jiayuan and Sch{\"o}lkopf, Bernhard},
  journal={Advances in Neural Information Processing Systems},
  volume={19},
  year={2006}
}

@article{klamt2009hypergraphs,
  title={Hypergraphs and cellular networks},
  author={Klamt, Steffen and Haus, Utz-Uwe and Theis, Fabian},
  journal={PLoS Computational Biology},
  volume={5},
  number={5},
  pages={e1000385},
  year={2009},
  publisher={Public Library of Science San Francisco, USA}
}

@inproceedings{feng2019hypergraph,
  title={Hypergraph neural networks},
  author={Feng, Yifan and You, Haoxuan and Zhang, Zizhao and Ji, Rongrong and Gao, Yue},
  booktitle={AAAI Conference on Artificial Intelligence},
  volume={33}, 
  pages={3558--3565},
  year={2019}
}

@article{gao2020hypergraph,
  title={Hypergraph learning: Methods and practices},
  author={Gao, Yue and Zhang, Zizhao and Lin, Haojie and Zhao, Xibin and Du, Shaoyi and Zou, Changqing},
  journal={IEEE Transactions on Pattern Analysis and Machine Intelligence},
  year={2020},
  publisher={IEEE}
}

@article{mukherjee2022estimation,
  title={Estimation in tensor Ising models},
  author={Mukherjee, Somabha and Son, Jaesung and Bhattacharya, Bhaswar B},
  journal={Information and Inference: A Journal of the IMA},
  volume={11},
  number={4},
  pages={1457--1500},
  year={2022},
  publisher={Oxford University Press}
}

\end{document}